\pdfoutput=1
\documentclass[11pt]{article}

\usepackage[letterpaper,margin=1in]{geometry}
\usepackage[T1]{fontenc}
\usepackage{amsmath,amssymb,amsfonts,amsthm}
\usepackage{newtxtext}
\usepackage{newtxmath}
\usepackage{mathtools,bm}
\usepackage{setspace}
\usepackage[round,comma]{natbib}
\bibpunct[, ]{(}{)}{,}{a}{}{,}

\usepackage{booktabs,array}
\usepackage{algorithm}
\usepackage{algpseudocode}
\usepackage{xcolor}
\usepackage{tikz}
\usepackage{pgfplots}
\pgfplotsset{compat=1.18}
\usepgfplotslibrary{groupplots}
\usetikzlibrary{positioning,shapes.geometric,arrows.meta,patterns}
\usepackage{hyperref}
\hypersetup{colorlinks=true,linkcolor=blue,citecolor=blue,urlcolor=blue,filecolor=blue,
 pdftitle={Constrained Nested-Logit Pricing: An Exact Reformulation and Global Optimization Method},
 pdfauthor={Hoang Giang Pham, Tien Mai}}

\theoremstyle{plain}
\newtheorem{theorem}{Theorem}
\newtheorem{proposition}{Proposition}
\newtheorem{corollary}{Corollary}

\newcommand{\R}{\mathbb{R}}
\newcommand{\cM}{\mathcal{M}}
\newcommand{\cN}{\mathcal{N}}
\newcommand{\bp}{\bm{p}}
\newcommand{\bu}{\bm{u}}

\newcommand{\bd}{\bm{d}}

\newcommand{\Kexp}{\mathcal K_{\mathrm{exp}}}

\title{Constrained Nested-Logit Pricing: An Exact Reformulation and Global Optimization Method}
\author{%
 Hoang Giang Pham \qquad Tien Mai\\[2pt]
 {\small School of Computing and Information Systems}\\
 {\small Singapore Management University}
}
\date{}

\begin{document}
\maketitle

\begin{abstract}
Sellers must price substitutable products while meeting operational requirements
such as price limits and sales targets. Nested logit captures stronger substitution among
similar products, but these requirements can undermine the structure used to solve
unconstrained pricing problems. We develop an exact reformulation and a global optimization
method for a fixed assortment with price bounds and weighted price caps. The key novelty is
to concentrate the difficult interactions at the nest level while retaining convex
optimization over individual products, without discretizing prices or approximating demand.
This structure makes the number of nests, rather than products, govern the exponent in a
conditional worst-case complexity bound. The method returns feasible prices together with a
bound on how much profit could still be gained, allowing decision makers to assess solution
quality even when computation stops early. Extensions cover purchase-probability bounds and
price-order rules, with explicit accounts of their computational cost. In synthetic
experiments, the method certifies solutions for assortments of up to $500$ products with
few nests. Projecting unconstrained prices onto the feasible set loses up to $32.60\%$ of profit
relative to direct constrained optimization, highlighting the value of jointly optimizing
substitution and feasibility.

\end{abstract}

\section{Introduction}

Pricing a portfolio of substitutable products requires more than balancing each product's
margin against its own demand. Raising one price can redirect customers to another product,
change the sales mix, or induce no purchase. A seller must therefore optimize prices jointly,
accounting for how customers substitute across alternatives. This decision is central to
revenue management in transportation, hospitality, retailing, and online services
\citep{gallego1994,talluri2004}.

The demand model determines which substitutions the pricing decision anticipates. The
multinomial logit (MNL) model is attractive computationally because its choice probabilities
have a closed form, but its independence-of-irrelevant-alternatives property imposes
proportional substitution across competing alternatives \citep[Chapter~3]{train2009}. This can be restrictive when
products form natural groups: two fares in the same family or two rooms of a similar class
may be closer substitutes than products in different groups. Nested logit (NL) represents
such groups as nests and allows correlated unobserved utilities within a nest while retaining
closed-form choice probabilities \citep{mcfadden1978,benakiva1985,train2009}. It offers a
practical way to account for differentiated substitution without requiring simulation-based
demand evaluation. That distinction matters for decisions: in our matched experiments,
prices optimized under a calibrated MNL approximation lose up to $34.42\%$ of profit when
evaluated under NL demand.

Operational constraints introduce a second requirement. Firms may impose price floors and
ceilings, cap a weighted price index, maintain a price ladder, or require a minimum conversion
rate \citep{davis2017,shaomai2026}. These restrictions couple products or tie pricing to aggregate demand. They can also
eliminate the common-markup (CM) or adjusted-markup structure of unconstrained NL pricing
\citep{lihuh2011,gallegowang2014}. CM pricing sets $p_j-c_j=t$ for all products;
adjusted-markup characterizations instead equalize a quantity that accounts for demand
sensitivities and nest parameters. For example, a binding ceiling can force one product's
markup below the common value while other products remain unconstrained. Different bounds can bind within the same
nest, while a purchase-probability target couples all nests through the outside option.
Computing unconstrained prices and then restoring feasibility need not preserve profit:
projected CM prices lose up to $32.60\%$ in our experiments. The operational
problem is therefore to optimize substitution and feasibility jointly.

A third requirement is to assess the quality of the prices returned. NL profit is generally
nonconcave in prices \citep{hansonmartin1996,gallegowang2014}, so a local optimizer can produce a feasible decision without bounding
how much better another feasible decision could be. A \emph{global optimality certificate}
combines a feasible price vector with a valid bound on the profit improvement still possible
\citep{tawarmalani2002}.
Such a bound supports accuracy-based stopping and lets a decision maker assess the remaining
optimization error even when a computation ends at a resource limit. It certifies the pricing
optimization problem for the specified demand model and constraints; it does not certify the
accuracy of the demand model itself.
Existing work provides markup characterizations for unconstrained NL pricing and global
approximation guarantees for several constrained variants
\citep{lihuh2011,gallegowang2014,rayfield2015,phammai2026}. The remaining gap addressed here is
an exact continuous formulation for coupled operational constraints whose global search
dimension is explicit in the number of nests. Section~\ref{sec:literature} compares the
constraint classes and computational approaches in detail.

Writing $N$ for the number of nests, our approach separates product-level convex optimization
from nest-level nonconvex search through aggregation within nests. Adding products enlarges the convex node model without
increasing the base number of bilinear links. In particular, the convex subproblems still contain product-level variables and grow
with product count, while the number of bilinear relations requiring spatial relaxation
depends on nest count. Here a bilinear link is a constraint involving multiplication of two
decision variables; it does not refer to a pair of products offered for sale. This offers a way to retain the computational
benefit of nest structure when operational constraints invalidate the usual markup formulas. Our paper makes four contributions.
\nopagebreak[4]
\begin{enumerate}
 \item \textbf{An exact reformulation with only $N+1$ nonconvex links.}
 We handle continuous prices, product-specific sensitivities, finite price intervals, and
 nonnegative linear price couplings without discretizing prices or approximating demand.
 One bilinear link represents the objective ratio and one represents each nest; the remaining
 nonlinear constraints have exact exponential-cone representations. \textit{To the best of our
 knowledge, this is the first exact reformulation for this constrained NL pricing class} to
 combine this link count with an explicit conditional fixed-nest approximation bound
 (Section~\ref{sec:reformulation}).

 \item \textbf{Certified pricing decisions and polynomial product-count dependence for fixed nests.}
 We develop spatial branch-and-bound (B\&B) with convex envelopes, bound propagation, and
 a persistent conic model. It returns feasible prices and a bound on remaining profit
 improvement, so solution quality can be assessed even after an early stop. Under the assumptions stated in
 Theorem~\ref{thm:bb-complexity}, balanced refinement, which systematically bisects both
 factor intervals of each bilinear link into the same number of pieces per factor, needs
 $O((C_{\rm B}/\varepsilon)^{N+1})$ nodes for absolute profit accuracy $\varepsilon>0$,
 with polynomial arithmetic work per node. Here $C_{\rm B}$ depends on the instance and
 initial bounds. For fixed $N$
 and polynomially controlled constants, total arithmetic work is polynomial in product
 count, constraint count, and inverse accuracy. A controlled positive feasible-profit scale
 also yields a relative approximation scheme. This provides a conditional polynomial-time
 approximation scheme (PTAS) interpretation in the arithmetic model
 (Theorem~\ref{thm:bb-complexity} and Corollary~\ref{cor:fixed-nest-approx}).
 These guarantees assume a uniform feasible-solution recovery property and controlled
 numerical conditioning; they are not established for every instance in the model class.
 Here an approximation scheme returns a feasible solution within a specified relative
 error of optimal profit; the fixed-nest result bounds arithmetic work polynomially in
 problem dimensions and inverse accuracy under the stated assumptions.

 \item \textbf{A precise account of the cost of operational constraints.}
 We identify which extensions preserve the small search dimension. Lower bounds on total
 purchase probability and price orders between products with equal sensitivities add no
 bilinear links. Exact upper purchase bounds add nest-level links, while detailed
 product-share requirements can make the link count grow with product count. This
 distinguishes restrictions that retain the base computational advantage from those that
 can make global certification more expensive (Section~\ref{sec:additional-constraints}).

 \item \textbf{Computational evidence and economically material pricing improvements.}
 B\&B certifies all $160$ matched instances, $69$ of $96$ larger instances, and every
 product-count case up to $500$ products with two or four nests. The remaining large
 instances have relative shifted-loss gaps at most $0.282\%$. The economic comparisons
 show why both nested substitution and direct constraint handling matter: MNL-based prices
 lose up to $34.42\%$ of profit under NL demand, and projected unconstrained CM
 prices lose up to $32.60\%$. These findings demonstrate the value of the proposed approach
 on the generated instances and document the growing cost of certification as nest count
 increases (Section~\ref{sec:experiments}).
\end{enumerate}

The paper is organized as follows. 
Section~\ref{sec:literature} reviews related work, Section~\ref{sec:model} defines the model,
Section~\ref{sec:reformulation} derives the reformulation, Section~\ref{sec:algorithm} presents
the algorithm and its complexity, Section~\ref{sec:additional-constraints} develops operational
extensions, Section~\ref{sec:experiments} reports the computational study, and
Section~\ref{sec:conclusion} concludes. The appendices contain the proofs, the complete B\&B
pseudocode, detailed experimental protocols, and supplementary computational results.

\section{Literature review}\label{sec:literature}

The review distinguishes results for assortment selection, unconstrained pricing, and
operationally constrained pricing, then identifies the global optimization tools used here.

\paragraph{Choice models and assortment optimization.}
NL belongs to the generalized extreme-value family of random-utility models
\citep{mcfadden1978,koppelmanwen1998,train2009}. In choice-based revenue management,
assortment optimization under NL admits polynomial algorithms for important nest-separable
cases and exact methods for more general variants
\citep{talluri2004,davis2014,gallegotopaloglu2014,alfandari2021}.
These results exploit nest structure, but concern which products to offer. Here the assortment
is fixed and the decisions are continuous, product-specific prices coupled by operational
restrictions, so assortment tractability does not directly resolve the pricing problem.

\paragraph{Unconstrained logit pricing.}
MNL profit is generally not concave in prices \citep{hansonmartin1996}, but market-share and
attraction transformations expose concavity or markup structure. For NL,
\citet{lihuh2011} establish concavity in market shares and characterize optimal prices under
common sensitivity conditions, and \citet{gallegowang2014} allow product-differentiated
sensitivities and arbitrary nest coefficients, derive adjusted-markup invariance, and reduce the
problem to scalar optimization. These low-dimensional characterizations provide the
unconstrained benchmark used in our experiments. Binding price or demand constraints can
destroy the markup structure on which they rely.

\paragraph{Constrained logit pricing.}
For NL with product-specific price intervals, \citet{rayfield2015} study
fixed-assortment and joint assortment--pricing problems, show that bounds can destroy scalar
unimodality, and develop a fully polynomial-time approximation scheme based on linear
programs, with computable upper bounds. \citet{chenjiang2020} develop approximation
methods for capacitated joint assortment and pricing, where capacity restricts the offered
assortment. These restrictions differ from linear couplings among continuous product prices.
\citet{davis2017} study quality-consistency constraints with prices
from a finite set, where the ordering structure yields polynomial-size linear programs. For MNL,
\citet{shaomai2026} handle price and purchase-probability constraints in static and network
settings with a bisection and mixed-integer programming method that controls a piecewise-linear
approximation error, and
\citet{phammai2026} develop branch-and-cut and piecewise-linear methods for constrained
assortment and pricing under the broader generalized nested logit model. The generalized model includes standard NL as a special case, so this work already provides
globally controlled approximations for a broader demand model. Our distinction is an exact
continuous formulation for the stated standard-NL constraint class, with the remaining
nonconvexity and the conditional approximation bound expressed through nest-level links.
The assumptions and arithmetic interpretation of that bound are stated in
Theorem~\ref{thm:bb-complexity} and Corollary~\ref{cor:fixed-nest-approx}.

\paragraph{Global optimization.}
Spatial B\&B relaxes nonconvex factors over boxes and recursively tightens them. McCormick
inequalities give the convex and concave envelopes of a bilinear term \citep{mccormick1976},
and \citet{tawarmalani2002} develop systematic algorithms for factorable nonconvex programs.
Exponential-cone representations allow the resulting convex node problems to be solved by
primal--dual interior-point methods \citep{nesterovnemirovskii1994,dahlandersen2022}.
These established tools support our algorithm; the paper-specific construction is the exact
NL reformulation that makes them applicable with a nest-dependent number of nonconvex links.

\section{Nested-logit demand and the  pricing problem}\label{sec:model}

We consider a single seller choosing prices for a fixed assortment in one selling period.
Demand parameters and costs are known, and the objective is expected profit; there is no
parameter learning or inventory evolution in the model. Customers choose one product or the
outside option after observing prices. We use the standard two-level NL choice model
\citep{mcfadden1978,train2009}, with the parameterization specified below.

Throughout, $\R$ denotes the real numbers, $\R_+$ the nonnegative real numbers, and
$\R_{++}$ the positive real numbers. Bold symbols collect scalar components into vectors,
and vector inequalities are componentwise. Let $\cM=[m]=\{1,\ldots,m\}$ be the set of products, where $m$ is the number of
products. The products are partitioned into nests
$\cN=\{1,\ldots,N\}$, with $\emptyset\ne\cM_n\subseteq\cM$ denoting the products in nest $n$, so that
$\cM=\bigcup_{n\in\cN}\cM_n$ and $\cM_n\cap\cM_{n'}=\emptyset$ for $n\ne n'$.
Let $n(j)$ be the nest containing product $j$. Product $j$ has price $p_j$, unit variable cost
$c_j$, and attraction intercept $a_j$. Let $\beta_j>0$ denote the product-specific price
sensitivity in the within-nest attraction. The within-nest attraction of product $j$ is $\exp(a_j-\beta_jp_j)$.
Each nest $n$ has a dissimilarity parameter $\sigma_n\in(0,1]$ and attraction sum
$W_n(\bp)=\sum_{j\in\cM_n}\exp(a_j-\beta_jp_j)$. Its logsum $\log W_n(\bp)$
is the inclusive value in the usual NL terminology \citep[Section~4.2]{train2009}.
We use the within-nest utility parameterization: in the conventional NL notation with
deterministic utility $V_j$, $a_j-\beta_jp_j=V_j/\sigma_{n(j)}$. Thus, $\beta_j$
already includes this scale adjustment; holding $\beta_j$ fixed while varying $\sigma_n$
holds within-nest attraction sensitivity fixed. This distinction between within-nest and
upper-level utility scaling is discussed by \citet{koppelmanwen1998} and
\citet[Section~4.2]{train2009}.
Equivalently, conventional utility is
$V_j=\sigma_{n(j)}a_j-\sigma_{n(j)}\beta_jp_j$.
Thus, varying $\sigma_n$ at fixed $\beta_j$ also changes the price coefficient in $V_j$;
it is not a comparison at fixed conventional utility coefficients.
The outside option is normalized to utility zero and forms a singleton alternative. The choice
denominator is $D(\bp)=1+\sum_{n\in\cN}W_n(\bp)^{\sigma_n}$.
Thus, the unconditional probability of choosing an alternative in nest $n$ and, for $j\in\cM_n$, the probability
of choosing product $j$ conditional on nest $n$ are $P_n^{\rm nest}(\bp)=W_n(\bp)^{\sigma_n}/D(\bp)$ and
$P_{j\mid n}(\bp)=\exp(a_j-\beta_jp_j)/W_n(\bp)$.
The unconditional purchase probability of product $j\in\cM_n$ is therefore
\begin{equation}
 P_j(\bp)=
 \frac{
 \exp(a_j-\beta_jp_j)W_n(\bp)^{\sigma_n-1}
 }{D(\bp)},
 \label{eq:product-probability}
\end{equation}
and the no-purchase probability is $P_0(\bp)=1/D(\bp)$. Notice that
$P_0(\bp)+\sum_{j\in\cM}P_j(\bp)=1$.

Parameter $\sigma_n$ controls substitution within nest $n$. When $\sigma_n$ is small, products in
that nest have strongly correlated unobserved utilities. If $\sigma_n=1$ for every nest, then
\eqref{eq:product-probability} reduces to the MNL model. For consistency with random-utility
maximization, the nesting structure and dissimilarity parameters must satisfy the usual
NL regularity conditions; the two-level specification above uses
$0<\sigma_n\le1$ \citep[Section~4.2.2]{train2009}.
Behaviorally, the nest represents products that customers regard as closer substitutes,
so a price increase can redirect demand toward other products in the same nest.
Comparisons across values of $\sigma_n$ must also specify which utility coefficients are
held fixed, as in the parameterization above.

We optimize expected profit per potential customer, denoted by
$\Pi(\bp)=\sum_{j\in\cM}(p_j-c_j)P_j(\bp)$. Multiplying this objective by a fixed
positive customer volume does not change the optimal prices, so no arrival-rate multiplier
is needed. Fixed costs are also omitted because they do not affect the decision. Substituting
\eqref{eq:product-probability} gives
\begin{equation}
 \Pi(\bp)=
 \frac{
 \displaystyle\sum_{n\in\cN}W_n(\bp)^{\sigma_n-1}
 \sum_{j\in\cM_n}(p_j-c_j)
 \exp(a_j-\beta_jp_j)
 }{
 \displaystyle 1+\sum_{n\in\cN}W_n(\bp)^{\sigma_n}
 }.
 \label{eq:fractional-profit}
\end{equation}
The static NL pricing problem is
\begin{equation}
 \max_{\bp}\bigl\{\Pi(\bp):\underline\bp\le\bp\le\overline\bp,\ B\bp\le\bd\bigr\}.
 \label{prob:nlp}
\end{equation}
Here $\underline p_j,\overline p_j\in\R$ are finite price bounds with
$\underline p_j\le\overline p_j$. The integer $n_B\ge0$ counts the coupling constraints,
$B\in\R_+^{n_B\times m}$ is their coefficient matrix, and $\bd\in\R^{n_B}$ is their
right-hand-side vector. Thus, each
row is a nonnegative budget, index, weighted-average, or related upper limit. For example,
$\sum_j\omega_jp_j\le d$ caps a price index with known weights $\omega_j\ge0$ and
known cap $d\in\R$; a weighted-average cap
uses weights summing to one. Explicit price intervals are kept separate from these coupling
rows because a lower price bound has a negative coefficient when written as a linear row.

Section~\ref{sec:additional-constraints} extends the model to total, nest-level, and
product-level purchase-probability requirements, as well as price orders and price-difference
rules. A total conversion target constrains $\sum_jP_j(\bp)$; a product-share requirement
constrains an individual $P_j(\bp)$. Price rules such as $p_i\le p_j$ contain mixed signs
and are treated separately: they become linear after transformation when the two products
have the same price sensitivity. The additional constraints can therefore be incorporated,
but their effect on the number of nonconvex links depends on their form.

For the base problem \eqref{prob:nlp}, we assume that
$\mathcal P=\{\bp\in\R^m:\underline\bp\le\bp\le\overline\bp,B\bp\le\bd\}$
is nonempty. The finite intervals make $\mathcal P$ compact, and continuity of the
choice probabilities guarantees that an optimal price vector exists. If nonnegative margins are
required, one can set $\underline p_j\ge c_j$.

The nonnegative-row condition on $B$ defines the base solver scope: these rows remain convex
after the change of variables in Section~\ref{sec:reformulation}. Generic signed rows and
coupling equalities need not have that property. Their exact logarithmic expressions remain
valid, but they can introduce nonconvexities beyond the nest-level links analyzed here.
The probability and common-sensitivity pairwise extensions are developed separately, with
explicit link counts in Section~\ref{sec:constraint-complexity}.

\section{Exact extended-variable reformulation}\label{sec:reformulation}

We reformulate the base pricing problem exactly, isolating its nonconvexity in $N+1$
bilinear links and expressing the remaining nonlinear constraints with exponential cones.
This structure enables the global solution method in Section~\ref{sec:algorithm}; operational
extensions are developed in Section~\ref{sec:additional-constraints}.

\subsection{A positive shifted-loss ratio}\label{sec:shifted-ratio}

We first make the ratio suitable for one-sided convex bounds. Direct profit terms can have
either sign when prices below cost are permitted; a sufficiently large shift makes every
loss contribution positive without changing the optimal prices. Choose a constant $\alpha$ satisfying
\begin{equation}
 \alpha>\max\left\{0,\max_{\bp\in\mathcal P}\max_{j\in\cM}(p_j-c_j)\right\}.
 \label{eq:alpha}
\end{equation}
For example, any $\alpha>\max\{0,\max_j(\overline p_j-c_j)\}$ is sufficient and
requires no preliminary nonlinear optimization. Define the shifted margin $\delta_j=\alpha-(p_j-c_j)>0$.
Define the shifted loss $L(\bp)=\alpha-\Pi(\bp)$. Since $\alpha$ is constant,
minimizing $L$ is equivalent to maximizing profit per customer. Moreover, profit is a
probability-weighted average of product margins and the zero outside-option margin, so
\eqref{eq:alpha} implies $L(\bp)>0$.
The term ``shifted loss'' names this mathematical objective; it is not an observed business
loss or regret. Positivity will let us overestimate the shifted numerator and underestimate
the denominator without understating their ratio. It will also allow multiplication of
positive lower bounds on the nest factors without reversing their direction.
The outside option is important in this identity. From \eqref{eq:fractional-profit},
$L(\bp)=\widetilde N(\bp)/D(\bp)$, where
\begin{equation}
 \widetilde N(\bp)=\alpha+
 \sum_{n\in\cN}W_n(\bp)^{\sigma_n-1}
 \sum_{j\in\cM_n}\delta_j
 \exp(a_j-\beta_jp_j).
 \label{eq:shifted-numerator}
\end{equation}
The first term, $\alpha$, is the shifted loss associated with no purchase. All terms in
\eqref{eq:shifted-numerator} are nonnegative.

\subsection{Exponential variables and nest aggregates}\label{sec:exponential-variables}

The shifted ratio is still expressed in prices through exponentials. We next use a bijective
change of variables that makes each attraction linear and each shifted-margin contribution
an entropy term of the form $u\log u$. Nest aggregation then separates product-level convex
terms from the interactions that require global optimization. For every $j\in\cM$, define
\begin{equation}
 u_j=\exp(\beta_j\delta_j)
 =\exp\left(\beta_j(\alpha+c_j-p_j)\right).
 \label{eq:u-definition}
\end{equation}
The inverse transformation recovers prices exactly:
\begin{equation}
 p_j=c_j+\alpha-\frac{1}{\beta_j}\log u_j.
 \label{eq:recover-price}
\end{equation}
The transformation is one-to-one, and larger $u_j$ corresponds to lower price $p_j$.
Its two goals are to make attraction linear, $\exp(a_j-\beta_jp_j)=q_ju_j$, and to express
the shifted-margin contribution as a positive multiple of the convex term $u_j\log u_j$.
Before constructing the transformed model, we tighten each explicit price interval by solving
the two linear programs $p_j^{\min}=\min_{\bp\in\mathcal P}p_j$ and
$p_j^{\max}=\max_{\bp\in\mathcal P}p_j$. These bounds are finite because $\mathcal P$ is
compact. The corresponding bounds on $u_j$ are
\begin{equation}
 \exp\left(\beta_j(\alpha+c_j-p_j^{\max})\right)
 \le u_j\le
 \exp\left(\beta_j(\alpha+c_j-p_j^{\min})\right).
 \label{eq:u-bounds}
\end{equation}
Introduce the positive constants $q_j=\exp\left(a_j-\beta_j(c_j+\alpha)\right)$ for $j\in\cM$.
It follows that $\exp(a_j-\beta_jp_j)=q_ju_j$. Thus, the demand attractions are affine
in the new variables. For each nest, define its total attraction and shifted-margin sum by
\begin{align}
 z_n&=\sum_{j\in\cM_n}q_ju_j, &&n\in\cN, \label{eq:z-def}\\
 y_n&=\sum_{j\in\cM_n}\frac{1}{\beta_j}
 q_ju_j\log u_j, &&n\in\cN. \label{eq:y-def}
\end{align}
Thus, $z_n=W_n(\bp)$ is the total attraction of nest $n$, while $y_n$ is its
shifted-margin-weighted exponential sum. The function $u\log u$ is convex on $u>0$
\citep[Section~3.1.5]{boydvandenberghe2004},
so each $y_n$ is a separable convex function of the product variables in its nest.
Substituting these aggregates into the shifted loss gives
$L(\bu)=(\alpha+\sum_{n\in\cN}z_n^{\sigma_n-1}y_n)/(1+\sum_{n\in\cN}z_n^{\sigma_n})$.

\subsection{Isolating the nonconvex bilinear links}\label{sec:auxiliary-model}

The purpose of the auxiliary model is to retain every feasible price vector while separating
the nonconvex interactions from the convex entropy and power relations. We replace the ratio
by an objective variable $r$ and use one-sided inequalities whose directions are justified by
the positivity established in Section~\ref{sec:shifted-ratio}. The positive
variables $h$ and $k$ bound the shifted numerator and denominator, respectively. For each
nest, $z_n$ retains its definition in \eqref{eq:z-def}, while $y_n$ becomes an upper bound
on the entropy sum in \eqref{eq:y-def}. Introduce positive variables $\rho_n,w_n,\eta_n$
for the corresponding power and product terms. Thus, \eqref{eq:y-def} describes the exact
value associated with a price vector; \eqref{con:y} below is its epigraph representation.
Consider
\begin{subequations}\label{prob:nlbc}
\begin{align}
 \min\quad &r\quad\text{s.t.}\quad\eqref{eq:u-bounds},\ \eqref{eq:z-def},
 \tag{\theparentequation}\label{con:objective}\\*
 &h\ge\alpha+\sum_{n\in\cN}w_n,\qquad
 k\le1+\sum_{n\in\cN}\eta_n, \label{con:h}\\*
 &h\le rk,\qquad w_n\ge\rho_ny_n &&(n\in\cN), \label{con:rk}\\*
 &y_n\ge\sum_{j\in\cM_n}(q_j/\beta_j)u_j\log u_j &&(n\in\cN), \label{con:y}\\*
 &\rho_n\ge z_n^{\sigma_n-1},\qquad \eta_n\le z_n^{\sigma_n} &&(n\in\cN),
 \label{con:rho}\\*
 &\sum_{j\in\cM}B_{\ell j}\bigl(c_j+\alpha-\beta_j^{-1}\log u_j\bigr)\le d_\ell
 &&(\ell=1,\ldots,n_B). \label{con:general-price}
\end{align}
\end{subequations}
The optimization variables are $\bu\in\R_{++}^m$, $h,k,r>0$, and
$z_n,y_n,\rho_n,w_n,\eta_n>0$ for $n\in\cN$. In particular, \eqref{eq:alpha}
implies $u_j>1$, so the entropy sums and $y_n$ are positive.

The variables have distinct roles. Variable $u_j$ is the transformed image of price $p_j$;
$z_n$ and $y_n$ aggregate attraction and shifted margin within nest $n$;
$\rho_n$ and $\eta_n$ represent powers of $z_n$; $w_n$ represents the product
$\rho_ny_n$; and $h$ and $k$ represent the shifted numerator and choice denominator.
Every feasible price vector has a representation with all auxiliary inequalities tight.
Conversely, any feasible extended point satisfies $r\ge L(\bp)$ for its recovered prices.
To see why inequality bounds still give an exact reformulation, let $Y_n(\bu)$ denote the
exact entropy sum in \eqref{eq:y-def}. Since $y_n\ge Y_n(\bu)>0$ and
$\rho_n\ge z_n^{\sigma_n-1}>0$, their product satisfies
$\rho_ny_n\ge z_n^{\sigma_n-1}Y_n(\bu)$. Hence $h\ge\widetilde N(\bp)$, $0<k\le D(\bp)$, and
$r\ge h/k\ge\widetilde N(\bp)/D(\bp)=L(\bp)$. Every feasible price vector admits equality throughout this chain. Thus, the inequalities
cannot create a better objective than the original problem, and they can represent its optimum.
Consequently, an optimum can be chosen with $r=L(\bp)$, and its profit per customer is
$\alpha-r$. Auxiliary variables need not be tight at every feasible extended point.

The structure of \eqref{prob:nlbc} separates tractable nonlinearities from a small number of
nonconvex bilinear links. Constraints \eqref{eq:u-bounds}, \eqref{eq:z-def}, and \eqref{con:h} are affine. Constraint \eqref{con:y} is convex because
$u\mapsto u\log u$ is convex on $u>0$. To interpret these inequalities, recall two
basic convex-set constructions: the epigraph is the region above a function's
graph, and the hypograph is the region below it. An epigraph of a convex function and a
hypograph of a concave function are convex sets. The power-function curvature facts
are also standard \citep[Section~3.1.5]{boydvandenberghe2004}. Since $0<\sigma_n\le1$, the two power bounds in \eqref{con:rho} define the epigraph of the convex function
$z_n\mapsto z_n^{\sigma_n-1}$ and the hypograph of the concave function
$z_n\mapsto z_n^{\sigma_n}$, respectively. The $N+1$ nonconvex bilinear links are
$rk$ and $\rho_ny_n$ ($n\in\cN$), grouped in \eqref{con:rk}. Thus, its nonconvexity is isolated in a number of linking relations that depends
on the number of nests rather than the number of products. The transformed general price
constraints \eqref{con:general-price} are convex when $B_{\ell j}\ge0$ for all $j$ in a given
row, because they then contain a nonnegative weighted sum of $-\log u_j$; for an arbitrary
matrix $B$, however, these constraints may introduce additional nonconvexity beyond the
$N+1$ bilinear links.

\begin{proposition}\label{prop:reformulation-equivalence}
Under the compactness assumption on $\mathcal P$ and \eqref{eq:alpha}, formulations \eqref{prob:nlp} and
\eqref{prob:nlbc} have the same optimal prices. If $r^*$ is the optimal value of
\eqref{prob:nlbc}, the optimal profit per customer is $\alpha-r^*$.
\end{proposition}
The proof is given in Appendix~\ref{app:reformulation-proof}.

\subsection{A solver-ready exponential-cone representation}\label{sec:bilinear-exp}

The preceding model identifies the global search variables, but a solver also needs an
explicit representation of the remaining convex constraints. We now express the entropy,
power, and transformed price relations using a single standard cone. This changes neither
the feasible prices nor the optimum. We use standard exponential-cone modeling identities
\citep[Chapter~5]{mosekcookbook}; the application to the pricing reformulation is derived
below. The
\(N+1\) bilinear links \(rk\) and \(\rho_ny_n\), \(n\in\cN\), are retained
explicitly, whereas every remaining nonlinear term is placed in a standard log-domain convex
set. For implementation, we use the closed exponential cone
$\Kexp=\operatorname{cl}\{(x,y,z)\in\R\times\R_{+}\times\R_+:y\exp(x/y)\le z\}$,
where $\operatorname{cl}$ denotes closure. All uses in the exact model below have positive
second coordinates, so they use the displayed exponential inequality directly.
This cone is a standard convex set supported by conic solvers, not a piecewise approximation
to the exponential function. Closure includes boundary limits of the set; those limits are
not needed for the positive-coordinate identities used here. Each cone block below means
one constraint requiring a triple of expressions to belong to this set. We use ``log-domain''
to describe representations derived from logarithmic identities, without requiring the
solver to handle a nonconvex logarithmic equality.
After these substitutions, a B\&B node contains only affine constraints, standard convex cones,
and McCormick envelopes.

\paragraph{Entropy terms.}
To represent the convex product-level sums without approximation, for every product \(j\in\cM\), introduce \(e_j\in\R\) and replace \eqref{con:y} by
\begin{equation}
 (-e_j,u_j,1)\in\Kexp,\ \ j\in\cM,
 \qquad
 y_n\ge\sum_{j\in\cM_n}\frac{q_j}{\beta_j}e_j,\ \ n\in\cN.
 \label{con:bec-entropy}
\end{equation}
The cone constraint is equivalent to \(e_j\ge u_j\log u_j\).
Indeed, its defining inequality is $u_j\exp(-e_j/u_j)\le1$. Taking logarithms gives
$\log u_j-e_j/u_j\le0$, which gives the stated entropy bound because $u_j>0$.

\paragraph{Power term \(\rho_n\ge z_n^{\sigma_n-1}\).}
Taking logarithms gives $\log\rho_n+(1-\sigma_n)\log z_n\ge0$. Introduce real logarithmic underestimators \(\psi_n^\rho\) and \(\psi_n^z\). The power
constraint is exactly represented, for every $n\in\cN$, by
\begin{subequations}\label{con:bec-rho}
\begin{align}
 &(\psi_n^\rho,1,\rho_n)\in\Kexp,\qquad (\psi_n^z,1,z_n)\in\Kexp,
 \label{con:bec-rho-a}\\
 &\psi_n^\rho+(1-\sigma_n)\psi_n^z\ge0. \label{con:bec-rho-c}
\end{align}
\end{subequations}
The first two constraints give \(\psi_n^\rho\le\log\rho_n\) and
\(\psi_n^z\le\log z_n\); hence the system is feasible if and only if the logarithmic
inequality holds.

\paragraph{Power term \(\eta_n\le z_n^{\sigma_n}\).}
Taking logarithms gives \(\log\eta_n\le\sigma_n\log z_n\). An exact
representation needs a different construction from the preceding paragraph: replacing
$\log\eta_n$ on the left by a lower bound could admit points that violate the intended
inequality. Instead, a solver-ready log-domain representation follows from the equivalent weighted-geometric-mean
inequality \(\eta_n\le z_n^{\sigma_n}1^{1-\sigma_n}\). Introduce real variables
\(s_{n1}^\eta,s_{n2}^\eta\) and impose, for every $n\in\cN$,
\begin{subequations}\label{con:bec-eta}
\begin{align}
 &(s_{n1}^\eta,\eta_n,z_n)\in\Kexp,\qquad (s_{n2}^\eta,\eta_n,1)\in\Kexp,
 \label{con:bec-eta-a}\\
 &\sigma_ns_{n1}^\eta+(1-\sigma_n)s_{n2}^\eta\ge0. \label{con:bec-eta-c}
\end{align}
\end{subequations}
Indeed, the cone constraints imply $s_{n1}^\eta\le\eta_n\log(z_n/\eta_n)$ and
$s_{n2}^\eta\le\eta_n\log(1/\eta_n)$, so the last inequality is feasible if and only if
\(\sigma_n\log z_n-\log\eta_n\ge0\).

\paragraph{Transformed price constraints.}
The last step preserves the nonnegative linear price couplings. If \(B_{\ell j}\ge0\) for every \(j\) in row \(\ell\), introduce \(\psi_j^u\in\R\) and replace
\eqref{con:general-price} by
\begin{equation}
 (\psi_j^u,1,u_j)\in\Kexp,\ \ j\in\cM,
 \qquad
 \sum_{j\in\cM}B_{\ell j}(c_j+\alpha)-\sum_{j\in\cM}\frac{B_{\ell j}}{\beta_j}\psi_j^u\le d_\ell,
 \ \ \ell=1,\ldots,n_B.
 \label{con:bec-price}
\end{equation}
This representation is exact for nonnegative rows of \(B\). A row with mixed-sign
coefficients is generally nonconvex after the logarithmic transformation and must be retained
as an additional nonconvex constraint or treated by a separate relaxation.

Combining these representations gives the solver-ready program
\begin{equation}
 \min r\quad\text{s.t.}\quad
 \eqref{eq:u-bounds},\ \eqref{eq:z-def},\ \eqref{con:h},\ \eqref{con:rk},\ \eqref{con:bec-entropy},\
 \eqref{con:bec-rho},\ \eqref{con:bec-eta},\ \eqref{con:bec-price}.
 \label{prob:bec}
\end{equation}
The original positive domains from \eqref{prob:nlbc} are retained. In particular,
$u_j,z_n,\rho_n,\eta_n>0$, so all logarithms and the positive-coordinate cone
equivalences above are well defined. Formulation \eqref{prob:bec} contains \(N+1\) displayed
bilinear links, $rk$ and $\rho_ny_n$ ($n\in\cN$), grouped in \eqref{con:rk}. Every other nonlinear
constraint belongs to \(\Kexp\). Hence, for nonnegative rows of \(B\),
\eqref{prob:bec} is an exact solver-ready extended-variable reformulation of
\eqref{prob:nlbc}; no binary discretization is required. The count is an upper bound on
the links that actually require branching: for example, when $\sigma_n=1$, one can fix
$\rho_n=1$ without excluding an optimal price vector, and its nest product becomes affine.

\section{Global solution method}\label{sec:algorithm}

This section develops a spatial B\&B method that returns feasible prices and a global bound
on remaining profit improvement. The method exploits the reformulation in
Section~\ref{sec:reformulation}: it relaxes the $N+1$ bilinear links to obtain convex node
problems, recovers feasible pricing decisions, and tightens the bounds through propagation
and branching. We present the node operations, persistent-model implementation, stopping
certificate, and conditional complexity analysis.

\subsection{Node relaxation, incumbent recovery, and branching}

Here a box is a collection of lower and upper intervals for the factors in the bilinear
links. At each node, we solve a convex relaxation over one box to obtain a lower bound on
shifted loss. Feasible prices give an upper bound; the best such prices found so far are
called the incumbent. If a box cannot improve that incumbent beyond the tolerance, we
discard it; otherwise we can split a factor interval and examine the two resulting boxes.
Bound propagation tightens intervals using implications of the constraints. Local polishing
means running a local optimizer from recovered prices to improve their objective or restore
feasibility. Best-bound order processes the open box with the smallest lower bound first.
Each node restricts the product factors to a box and bounds the best shifted loss attainable
there. Initial bounds follow from the compact price system. During propagation, we restrict
attention to exact representations of price vectors, preserving at least one representation
of every feasible price. We use the identities $\rho_n=z_n^{\sigma_n-1}$,
$\eta_n=z_n^{\sigma_n}$, and $k=1+\sum_n\eta_n$, together with attraction-weighted bounds
linking $y_n$ and $z_n$; the incumbent also bounds $r$ from above.

At a B\&B node, let $v_0$ represent $rk$ and let $v_n$ represent
$\rho_ny_n$ for $n\in\cN$. For any one of these links, write its factors as $x$ and $y$,
with current bounds $[\underline x,\overline x]$ and
$[\underline y,\overline y]$. The relaxation of $v=xy$ uses the four convex-envelope
inequalities \citep{mccormick1976}
$v\ge\underline x y+\underline y x-\underline x\underline y$,
$v\ge\overline x y+\overline y x-\overline x\overline y$,
$v\le\overline x y+\underline y x-\overline x\underline y$, and
$v\le\underline x y+\overline y x-\underline x\overline y$.
In \eqref{con:rk}, replace $h\le rk$ by $h\le v_0$ and $w_n\ge\rho_ny_n$ by $w_n\ge v_n$. The four inequalities then link $v_0$ to $(r,k)$
and $v_n$ to $(\rho_n,y_n)$. A node is pruned when its lower bound rules out improvement
beyond the prescribed tolerance, with bound accounting as described below.
For the denominator product, only an overestimate $v_0-rk$ can favor the relaxation; for a
nest numerator product, only an underestimate $\rho_ny_n-v_n$ can do so. The branching score
therefore uses the positive part of the appropriate residual, normalized by
$\max\{1,|xy|\}$. If all realized residuals are small, the algorithm selects the link with
the largest worst-case envelope error rather than discarding the node. A candidate incumbent is
obtained by recovering prices from the node's $\bu$ solution and evaluating the original
fractional objective. It is accepted only if all original constraints hold; otherwise
local polishing must first find a feasible point. Local
polishing is invoked selectively, avoiding a full nonlinear solve at every node. The algorithm
terminates when the incumbent-to-global-bound gap, including tolerance-pruned boxes, meets
the prescribed absolute or relative tolerance.
We reuse one MOSEK Fusion model, updating factor bounds and envelope coefficients at each
node to save model-construction work \citep{mosekfusion}. This reuse does not require an
interior-point warm-start guarantee; implementation details are in Appendix~\ref{app:algorithm}.

\subsection{Complete procedure and certificate}

Appendix~\ref{app:algorithm} gives the complete procedure and bound bookkeeping.
Let $\underline L$ be the global lower bound on shifted loss and
$\overline L=L(\widehat\bp)$ the upper bound from the best feasible prices found.
We define the relative optimality gap as
\begin{equation}
 \operatorname{gap}(\underline L,\overline L)
 =\frac{\max\{0,\overline L-\underline L\}}
 {\max\{1,|\overline L|\}}.
 \label{eq:bb-relative-gap}
\end{equation}
If $\Pi^*$ denotes optimal profit per customer, these bounds imply
$\alpha-\overline L\le\Pi^*\le\alpha-\underline L$ and
$\Pi^*-\Pi(\widehat\bp)\le\overline L-\underline L$. Thus, the absolute shifted-loss gap equals a bound on the remaining profit improvement
per customer; no rescaling is required. The
relative gap in \eqref{eq:bb-relative-gap} is normalized by shifted loss, which depends on
$\alpha$; it is not a relative profit-loss bound. For example, if $\alpha=10$,
$\overline L=8$, and $\underline L=7.99$, the relative shifted-loss gap is $0.125\%$.
The incumbent profit is $2$, with at most $0.01$ additional profit possible, or $0.5\%$
of the incumbent profit. Thus, a shifted-loss percentage should not be read as a percentage
of profit. We use ``certified'' to mean that the
stated optimization tolerance has been met. Numerical certificates rely on the conic solver's
feasibility and dual tolerances; the implementation does not provide interval-arithmetic
verification.

The absolute and relative stopping tolerances are $\varepsilon_{\rm abs}>0$ and
$\varepsilon_{\rm rel}>0$. The global lower bound accounts for both open boxes and boxes
discarded within the absolute tolerance; exhausting the tree therefore need not prove zero
error. Appendix~\ref{app:algorithm} specifies this accounting. Incumbent updates require
original-model feasibility. The residual threshold $\varepsilon_{\rm bil}>0$ guides branching
only: small bilinear residuals indicate agreement with multiplication equations at one
relaxation point, whereas stopping requires a valid incumbent-to-global-bound gap.

\subsection{Worst-case arithmetic complexity}\label{sec:algorithm-complexity}

We analyze the base formulation \eqref{prob:bec}, with $d_{\rm L}=N+1$ bilinear links
and finite initial factor bounds. The analysis assumes the following
\emph{uniform relaxation-and-recovery (URR)} condition: any node whose maximum
McCormick-envelope error is at most $\theta>0$ can either be certified infeasible or yield
feasible prices with shifted loss at most $\kappa d_{\rm L}\theta$ above its relaxation
bound. The constant $\kappa>0$ is uniform over nodes and $\theta$, and this check or recovery
requires polynomial arithmetic work. The bound concerns balanced factor bisection, with
absolute shifted-loss tolerance $\varepsilon>0$.
In plain terms, URR assumes that sufficiently accurate bilinear envelopes lead either to a
proof of infeasibility or to feasible prices whose objective is comparably accurate. It is
an additional assumption, not a recovery guarantee proved here for every pricing instance.
The algorithm's selective local polishing does not by itself establish URR.
Let $n_v$ and $n_e$ denote the numbers of scalar variables and exponential-cone blocks in a
node relaxation. Write $n_K=O(n_v+n_B+n_e)$ for its Newton-system dimension and
$\vartheta=O(n_v+n_B+3n_e)$ for its barrier parameter. The exponential-cone barrier is
described by \citet{dahlandersen2022}; the iteration bound below uses standard
self-concordant-barrier theory \citep{nesterovnemirovskii1994}. We assume the interior-point regularity
and conditioning needed to solve each relaxation to accuracy
$\tau=\Theta(\varepsilon)$ in
$O(\sqrt{\vartheta}\log(\vartheta/\varepsilon)n_K^3)$ arithmetic operations with dense
linear algebra. Let $T_{\rm aux}(\varepsilon)$ be a polynomial bound on the additional
per-node work for propagation, recovery, and feasibility checks. Recovery and solve-error
budgets are chosen to total at most $\varepsilon$.
Here $n_K$ measures the size of the linear system solved in each interior-point iteration,
and the barrier parameter $\vartheta$ controls the iteration bound. The detailed expression
separates this convex-solve cost from the number of boxes searched globally.

\begin{samepage}
\begin{theorem}[Worst-case complexity of spatial B\&B]\label{thm:bb-complexity}
Under the conditions above, balanced spatial B\&B certifies an incumbent-to-global-bound
gap of at most $\varepsilon$ in shifted loss using
\begin{equation}
 O\!\left[\left(\frac{C_{\rm B}}{\varepsilon}\right)^{N+1}\right]
 \label{eq:bb-node-complexity}
\end{equation}
nodes and overall arithmetic work
\begin{equation}
 O\!\left[
 \left(\frac{C_{\rm B}}{\varepsilon}\right)^{N+1}
 \left\{\sqrt{\vartheta}\log\!\left(\frac{\vartheta}{\varepsilon}\right)n_K^3
       +T_{\rm aux}(\varepsilon)\right\}
 \right].
 \label{eq:bb-total-complexity}
\end{equation}
Here $C_{\rm B}$ depends on $\kappa$, $N$, and the initial factor widths, but not on
$\varepsilon$.
\end{theorem}
\end{samepage}
The proof is given in Appendix~\ref{app:complexity-proof}.
For the base model, $n_v$ and $n_e$ grow linearly with $m+N$, apart from auxiliaries for
$n_B$ price rows. Thus, for fixed $N$ and polynomially controlled constants and per-node
costs, the bound is polynomial in product count, constraint count, and inverse accuracy.
Its search exponent is $N+1$. If auxiliary work is dominated by the conic solve, the
$T_{\rm aux}$ term can be absorbed into the conic cost. This is a conditional arithmetic
result, not an unconditional bit-complexity guarantee.
The arithmetic model counts numerical operations without fully accounting for the number
of bits needed to represent their operands. ``Polynomially controlled'' means that the
instance-dependent constants must obey bounds that grow at most polynomially across the
instance family, rather than merely being finite for each instance. Conditioning describes
how sensitive the numerical solves and recovery are to small perturbations; poor conditioning
can invalidate the assumed uniform work bounds.
This is a worst-case bound for a balanced refinement schedule; the residual-based branching
rule in Algorithm~\ref{alg:spatial-bb} can give a different tree size. Bound tightening,
incumbent pruning, and sparsity can reduce observed work. Section~\ref{sec:additional-constraints}
counts the links added by operational extensions.

\begin{corollary}[Conditional fixed-nest relative approximation]\label{cor:fixed-nest-approx}
Fix $N$ and consider a family of base-model instances satisfying the conditions of
Theorem~\ref{thm:bb-complexity}. Suppose a feasible price vector $\bp^0$ with positive
profit $\Pi(\bp^0)>0$ is available in polynomial arithmetic work. Let the desired
relative tolerance be $\gamma\in(0,1)$. Suppose $C_{\rm B}/\Pi(\bp^0)$ is polynomially bounded
in the input dimensions and node-solve work at absolute accuracy proportional to $\gamma \Pi(\bp^0)$
is polynomial in those dimensions and $\log(1/\gamma)$. Then balanced refinement returns a feasible
$\widehat\bp$ satisfying $\Pi(\widehat\bp)\ge(1-\gamma)\Pi^*$ in arithmetic work polynomial in the numbers of products and constraints and in $1/\gamma$.
\end{corollary}
The proof is given in Appendix~\ref{app:complexity-proof}.
This corollary supplies a relative, rather than merely additive, approximation interpretation
of the fixed-nest result. Its assumptions must hold uniformly over the instance family:
finite constants for each separate instance are insufficient. Because the bound counts
arithmetic operations under the assumptions of Theorem~\ref{thm:bb-complexity}, it is a conditional
arithmetic approximation scheme, not an unconditional PTAS in binary input length.
Individual-share extensions can also increase the exponent, as described next.

\section{Additional operational constraints}\label{sec:additional-constraints}

Price bounds and weighted price caps alone do not capture every operational requirement:
a feasible pricing policy may also need to meet demand targets or maintain a coherent price
ladder across related products \citep{shaomai2026,davis2017}. Two further classes of restrictions
fit the extended-variable architecture: constraints on purchase probabilities and pairwise price constraints. Probability lower bounds and upper
bounds require different treatments because the base model only bounds the denominator from
above. Pairwise price rules become linear in the transformed variables when sensitivities are
common across products.

\subsection{Constraints on purchase probabilities}

Purchase-probability bounds let a seller control demand as well as margins. A lower bound on
total conversion can enforce a minimum expected sales target, while product- or nest-level
bounds can maintain demand for a designated product line. Upper bounds can represent limited
inventory or service capacity: multiplying choice probabilities by a known customer volume
gives expected sales, and weighting them by resource consumption gives expected resource use
\citep{shaomai2026}. In this static model, such bounds control expected demand; they do not
guarantee that realized stochastic demand stays below capacity. Aggregate rows can also
restrict the expected sales mix across product groups. For example, one may require
\begin{equation}
 \underline\pi_j\le P_j(\bp)\le\overline\pi_j,
 \qquad j\in\mathcal J\subseteq\cM,
 \label{eq:product-share-bounds}
\end{equation}
or, more generally,
\begin{equation}
 \underline b_s\le\sum_{j\in\cM}A_{sj}P_j(\bp)\le\overline b_s,
 \qquad s\in\mathcal S.
 \label{eq:linear-share-constraints}
\end{equation}
Here $\mathcal J\subseteq\cM$ indexes products with requested share bounds,
$0\le\underline\pi_j\le\overline\pi_j\le1$ are their limits, $\mathcal S$ is a finite
set of aggregate rows, $A_{sj}\in\R$ are known coefficients, and
$\underline b_s\le\overline b_s$ are known row limits. The coefficients may aggregate products, encode
expected resource consumption, or compare the expected sales of different product groups.
Let $\mathcal J$ contain exactly the products whose probabilities appear in the requested
constraints, and let $\mathcal T\subseteq\cN$ contain the nests whose shares are needed.
For a general row in \eqref{eq:linear-share-constraints}, $\mathcal J$ includes every
$j$ with $A_{sj}\ne0$. Introduce product-share variables $\pi_j^{\rm prod}\in[0,1]$ and unnormalized shares
$\zeta_j>0$ for $j\in\mathcal J$, and nest-share variables
$\pi_n^{\rm nest}\in[0,1]$ for $n\in\mathcal T$. The superscripts distinguish
product shares from nest shares even when their numeric indices coincide.
In addition to \eqref{prob:bec}, impose the applicable equations below:
\begin{subequations}\label{con:probability-linking}
\begin{align}
 k&=1+\sum_{n\in\cN}\eta_n, \label{con:probability-k}\\
 \eta_n&=\rho_nz_n,
 &&n\in\cN, \label{con:probability-power-link}\\
 k\pi_n^{\rm nest}&=\eta_n,
 &&n\in\mathcal T, \label{con:probability-nest}\\
 \zeta_j&=q_ju_j\rho_{n(j)},\qquad k\pi_j^{\rm prod}=\zeta_j,
 &&j\in\mathcal J. \label{con:probability-product}
\end{align}
\end{subequations}
Together with the corresponding log-domain convex power constraints
\eqref{con:bec-rho} and \eqref{con:bec-eta}, relation
\eqref{con:probability-power-link} forces
\(\rho_n=z_n^{\sigma_n-1}\) and \(\eta_n=z_n^{\sigma_n}\). Consequently,
\eqref{con:probability-nest} and \eqref{con:probability-product} give $\pi_n^{\rm nest}=P_n^{\rm nest}(\bp)$ and
$\pi_j^{\rm prod}=P_j(\bp)$. Constraints \eqref{eq:product-share-bounds} and
\eqref{eq:linear-share-constraints} are therefore linear in the new share variables.
The new nonconvexities are bilinear links: $N$ in
\eqref{con:probability-power-link}, $|\mathcal T|$ in \eqref{con:probability-nest},
and $2|\mathcal J|$ in
the two product-share equations in \eqref{con:probability-product}. Each admits the same
convex-envelope relaxation used in Section~\ref{sec:algorithm}.
Total purchase-probability restrictions require fewer variables. Write
$D=1+\sum_n z_n^{\sigma_n}$ for the true choice denominator. Then
$P_{\mathrm{buy}}(\bp)=\sum_{j\in\cM}P_j(\bp)=1-1/D$. For $0\le\underline\pi\le\overline\pi<1$, a lower bound
$P_{\mathrm{buy}}\ge\underline\pi$ can be imposed exactly by
$k\ge 1/(1-\underline\pi)$ in the base formulation: since $k\le D$, this implies
$D\ge1/(1-\underline\pi)$, and any qualifying price vector admits $k=D$.
An upper bound, or a two-sided band, requires $k=D$. Add
\eqref{con:probability-k}--\eqref{con:probability-power-link}; these equations and the
convex power constraints force $\eta_n=z_n^{\sigma_n}$ and $k=D$. The band is then
equivalent to
\begin{equation}
 \frac{1}{1-\underline\pi}\le k\le
 \frac{1}{1-\overline\pi}.
 \label{eq:total-share-k-bounds}
\end{equation}
The upper bound therefore adds $N$ bilinear links; the lower bound alone adds none.
The asymmetry comes from $k\le D$ in the base formulation: a lower bound on $k$ forces
the same lower bound on $D$, but an upper bound on $k$ does not force an upper bound on
$D$. The added equalities remove this freedom by making $k$ the true denominator.
Nest-level restrictions additionally use \eqref{con:probability-nest} only for the nests
whose shares are constrained. Product-level variables are unnecessary in either case.

\subsection{Pairwise price constraints}

Pairwise price rules preserve the intended quality ordering of a product portfolio. For
example, a seller may require a premium room or a more flexible fare to cost at least as much
as its basic counterpart. Such quality-consistency constraints are studied explicitly in NL
pricing by \citet{davis2017}. Bounds on price differences can additionally limit the surcharge
between adjacent tiers. These rules couple prices across products and need not be respected
by an unconstrained profit optimum. Suppose the seller requires a price ordering
\begin{equation}
 p_i\le p_j,
 \qquad (i,j)\in\mathcal E,
 \label{eq:pairwise-price-order}
\end{equation}
where \(\mathcal E\) is a specified collection of product pairs. Assume in this subsection
that \(\beta_j=\beta>0\) for every product. From
\eqref{eq:recover-price}, constraint \eqref{eq:pairwise-price-order} is equivalent to
$c_i+\alpha-\beta^{-1}\log u_i\le c_j+\alpha-\beta^{-1}\log u_j$, and hence to
\begin{equation}
 u_i\ge \exp\!\left(\beta(c_i-c_j)\right)u_j,
 \qquad (i,j)\in\mathcal E.
 \label{con:pairwise-u-linear}
\end{equation}
Thus, pairwise price ordering becomes a collection of linear constraints in \(\bu\). It adds
neither log-domain nonlinearities nor bilinear links to \eqref{prob:bec}, and it can be imposed
directly at every B\&B node. Price-difference rules such as
\(p_i-p_j\le\Delta_{ij}\) admit the same treatment, where $\Delta_{ij}\in\R$ is the prescribed maximum price difference, with the coefficient in
\eqref{con:pairwise-u-linear} replaced by
\(\exp(\beta(c_i-c_j-\Delta_{ij}))\).

When the sensitivities are heterogeneous, the corresponding condition is
$\beta_i^{-1}\log u_i-\beta_j^{-1}\log u_j\ge c_i-c_j$, which is generally nonconvex because it is a difference of concave logarithms. It must then
be retained as a separate nonconvex constraint or relaxed through bounded logarithmic
variables. The common-sensitivity assumption is therefore what makes pairwise price rules
especially inexpensive in the extended-variable formulation.

\subsection{Complexity of the extensions}\label{sec:constraint-complexity}

Theorem~\ref{thm:bb-complexity} applies to an extended model after replacing $d_{\rm L}=N+1$
by its total number of spatially relaxed bilinear links, provided the extended model has
finite factor bounds and satisfies the same technical conditions as the base model.
The conic node size grows polynomially with the added variables and constraints. A lower
bound on total purchase probability and common-sensitivity pairwise price rules add no
links, so their worst-case node-count exponent remains $N+1$. An upper bound on total
purchase probability adds the $N$ exact-power links
\eqref{con:probability-power-link}, giving $d_{\rm L}=2N+1$.

If exact shares are needed for nests in $\mathcal T$ and products in $\mathcal J$, the
full denominator and share representation \eqref{con:probability-linking} has
$d_{\rm L}=2N+1+|\mathcal T|+2|\mathcal J|$ bilinear links. Let $C_{\rm B}^{\rm ext}$ be the extension-specific analogue of $C_{\rm B}$. Under the stated
conditions, the worst-case node count is $O((C_{\rm B}^{\rm ext}/\varepsilon)^{d_{\rm L}})$,
with polynomial arithmetic work per node.
Thus, nest-share constraints preserve an exponent governed by $N$, whereas requesting
individual shares for a growing number of products can make the exponent grow with $m$.
The bound describes the displayed formulation; exploiting special aggregation structure
could reduce the required link count.

\section{Computational study}\label{sec:experiments}

The main text presents the experimental setup, evaluates certification across constraint
systems and demand models, and examines how product and nest counts affect computational
performance. The appendix provides supporting tables and detailed protocols
(Appendices~\ref{app:detailed-results} and~\ref{app:experimental-protocol}), together with
comparisons of NL and MNL pricing (Appendix~\ref{app:economic-effect}), a nest-scale ablation
(Appendix~\ref{sec:sigma-ablation}), the value of direct constrained optimization
(Appendix~\ref{sec:constant-markup-experiment}), and challenging instances for local optimization
(Appendix~\ref{sec:local-optimum-results}).

\subsection{Experimental design, benchmarks, and verification}\label{sec:experimental-design}

All experiments use synthetic instances, with profit measured per potential customer.
Tables report results to two decimal places; comparisons and certification use unrounded values.
We test linear price constraints (LC), total purchase-probability constraints (PPC), their
combination (LC+PPC), and pairwise price-order constraints (POC). Nests have approximately equal
sizes; their scales lie in strong-correlation $[0.25,0.50]$, medium $[0.50,0.75]$, or near-MNL
$[0.75,1]$ regimes. Intercepts are calibrated to a $0.50$ purchase probability at midpoint
prices $\bp^{\rm mid}=(\underline\bp+\overline\bp)/2$. Matched MNL instances set every
$\sigma_n=1$ and recalibrate intercepts to the same midpoint conversion.
We write $\mathrm{Uniform}[a,b]$ for a uniform distribution and $\mathrm{Normal}(\mu,s^2)$ for a
normal distribution with mean $\mu$ and variance $s^2$.
We compare B\&B with multistart sequential least squares programming (SLSQP), a local
constrained optimizer implemented in SciPy \citep{kraft1994,virtanen2020}, and, for $m\le10$, differential evolution (DE) with
local polishing \citep{stornprice1997}. These benchmarks supply feasible incumbents but no
global bound. Unless stated otherwise, B\&B uses relative and absolute tolerances $10^{-4}$
and $10^{-6}$, a $600$-second limit, and a $2500$-open-node limit. Open nodes are queued boxes;
solved nodes count processed relaxations. All reported gaps use the shifted-loss definition
\eqref{eq:bb-relative-gap}, rather than profit loss. For PPC, the implementation uses a weaker
outer relaxation and checks recovered prices against the original band; its bounds remain
valid, but the complexity result for the exact band formulation does not automatically apply.

Runs use MOSEK~11.2.4 \citep{mosekfusion} and SciPy~1.16.3 \citep{virtanen2020}
on an Intel Core i7-14700 workstation with 64~GB RAM.
The reported post-run audit checks incumbent feasibility, objectives, and gaps for all $356$
B\&B runs; it does not independently verify the conic dual bounds. Appendix~\ref{app:experimental-protocol}
gives the complete generator, campaign design, solver settings, and audit details;
Appendix~\ref{app:detailed-results} provides supporting numerical tables.

\subsection{Certification and demand-model fidelity}
\label{sec:final-results}

 This experiment measures how problem size, constraint system, and demand
model affect the ability of B\&B to bound solution quality. It also supplies the matched
instances for the economic comparison in Appendix~\ref{app:economic-effect}.
All $160$ matched runs, $80$ MNL and $80$ NL instances, reach the certified gap test,
and the median final gap is $0.00100\%$. Table~\ref{tab:final-family-performance} in
Appendix~\ref{app:detailed-results} reports median wall times and solved nodes by constraint
system, demand model, and instance size.
Demand-model structure helps explain the difference: at $\sigma_n=1$, the nest-power
relations simplify and the nest products can be made affine. MNL has median
runtime $0.20$ seconds and median two solved nodes, while NL has median runtime
$1.14$ seconds and median $160$ nodes. Within NL, nest count matters most: moving
from $N=2$ to $N=4$ raises median time from $0.52$ to $4.46$ seconds and median solved nodes
from $60$ to $696$, and the slowest run takes $119.5$ seconds. Stronger within-nest
correlation also enlarges the tree: strong-correlation instances need median $175$ solved
nodes, compared with $108$ for near-MNL instances. The most demanding cell is PPC with
$(m,N)=(20,4)$, at median $3922$ solved nodes.

\paragraph{B\&B performance across problem sizes and constraint systems.}

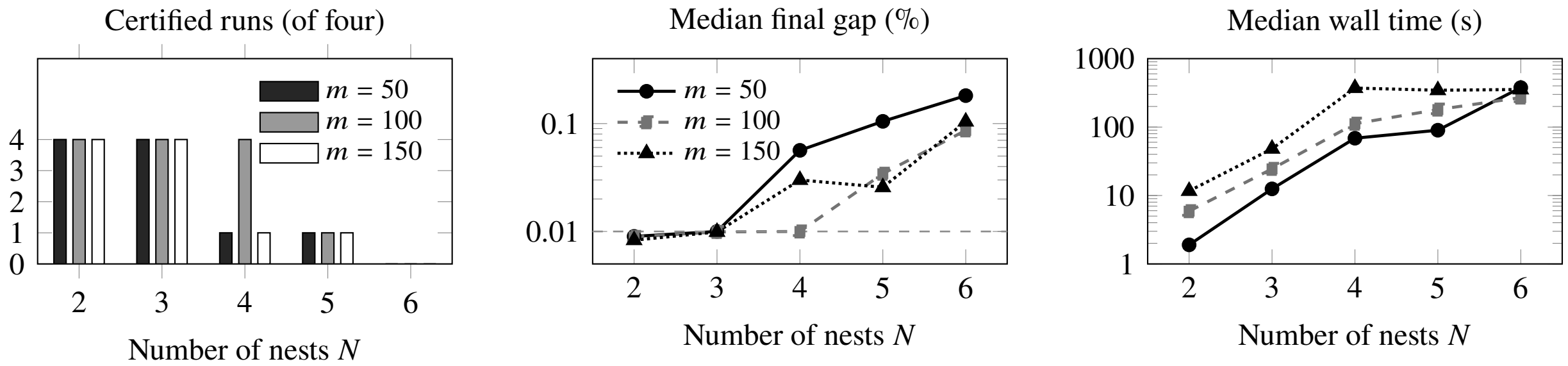
\begin{figure}[htbp]
\centering
\begin{tikzpicture}
\begin{groupplot}[
 group style={group size=3 by 1, horizontal sep=1.45cm},
 width=0.355\textwidth, height=3.7cm,
 xmin=1.5, xmax=6.5, xtick={2,...,6}, xlabel={Number of nests $N$},
 tick label style={font=\footnotesize}, label style={font=\footnotesize},
 title style={font=\footnotesize, yshift=-3pt},
 legend style={font=\scriptsize, draw=none, fill=none, row sep=-2pt},
 legend cell align=left]
\nextgroupplot[title={Certified runs (of four)}, ybar, bar width=3.6pt, ymin=0, ymax=6.6,
 ytick={0,...,4}, legend pos=north east, area legend]
\addplot[draw=black, fill=black!85] coordinates {(2,4) (3,4) (4,1) (5,1) (6,0)};
\addplot[draw=black, fill=black!40] coordinates {(2,4) (3,4) (4,4) (5,1) (6,0)};
\addplot[draw=black, fill=white] coordinates {(2,4) (3,4) (4,1) (5,1) (6,0)};
\legend{$m=50$, $m=100$, $m=150$}
\nextgroupplot[title={Median final gap (\%)}, ymode=log, ymin=0.005, ymax=0.4,
 ytick={0.01,0.1}, yticklabels={0.01,0.1}, legend pos=north west]
\addplot[black, mark=*, line width=0.9pt, mark size=1.7pt]
 coordinates {(2,0.0090) (3,0.0100) (4,0.0565) (5,0.1049) (6,0.1817)};
\addplot[black!55, dashed, mark=square*, line width=0.9pt, mark size=1.6pt]
 coordinates {(2,0.0089) (3,0.0099) (4,0.0100) (5,0.0343) (6,0.0865)};
\addplot[black, densely dotted, mark=triangle*, mark options={solid}, line width=0.9pt,
 mark size=2pt] coordinates {(2,0.0083) (3,0.0099) (4,0.0300) (5,0.0257) (6,0.1043)};
\draw[black!45, dashed, line width=0.5pt] (axis cs:1.5,0.01) -- (axis cs:6.5,0.01);
\legend{$m=50$, $m=100$, $m=150$}
\nextgroupplot[title={Median wall time (s)}, ymode=log, ymin=1, ymax=1000,
 ytick={1,10,100,1000}, yticklabels={1,10,100,1000}]
\addplot[black, mark=*, line width=0.9pt, mark size=1.7pt]
 coordinates {(2,1.90) (3,12.50) (4,68.75) (5,89.95) (6,377.25)};
\addplot[black!55, dashed, mark=square*, line width=0.9pt, mark size=1.6pt]
 coordinates {(2,5.90) (3,24.50) (4,112.08) (5,180.18) (6,269.97)};
\addplot[black, densely dotted, mark=triangle*, mark options={solid}, line width=0.9pt,
 mark size=2pt] coordinates {(2,11.64) (3,48.17) (4,371.46) (5,345.42) (6,354.43)};
\end{groupplot}
\end{tikzpicture}
\caption{NL runs testing B\&B performance across problem sizes and constraint systems. Each point summarizes the four
constraint systems of one $(m,N)$ cell; the dashed horizontal line marks the $0.01\%$ target.
Uncertified runs stop at the open-node limit. For $N\ge4$, their times and gaps are
observed at that stop, rather than at certification.}
\label{fig:large-instance-performance}
\end{figure}

To test B\&B certification performance as product and nest counts vary, we conduct $96$ runs
with $m\in\{50,100,150\}$. For NL,
$N\in\{2,\ldots,6\}$, all instances use the strong-correlation regime, and each $(m,N)$ cell
contains LC, PPC, LC+PPC, and POC. The calibrated MNL runs pool the twelve corresponding runs
at each $m$ from the original $N\in\{2,3,4\}$ design; MNL is not rerun for $N=5,6$ because
nest labels have no meaning in that model. The target remains $10^{-4}$ and the time limit is
$3600$ seconds. The open-node limit is $2500$, except that the $(m,N)=(50,6)$ cell uses
$10000$. Figure~\ref{fig:large-instance-performance} summarizes the NL runs.

B\&B certifies $69$ of these $96$ runs. All $36$ MNL cases certify, with median times
of $1.3$, $4.8$, and $8.6$ seconds for $m=50$, $100$, and $150$. For NL, all $24$
runs with $N\le3$ certify, as do six of twelve with $N=4$ and three of twelve with $N=5$; none
of the twelve $N=6$ runs certifies. The POC runs certify in every $N=4$ and $N=5$ cell, and
at $(m,N)=(100,4)$ all four constraint systems certify. Every uncertified run stops at the
open-node limit, well before the time limit, with final gap at most $0.282\%$; the largest
median gap is $0.182\%$ for $(m,N)=(50,6)$. Increasing $m$ raises node cost: at $N=2$, median
time grows from $1.9$ seconds for $m=50$ to $11.6$ seconds for $m=150$, while median node
counts remain near $80$. Increasing $N$, by contrast, expands the tree: at $m=100$, where every
run through $N=4$ certifies, the median solved-node count rises from $76$ for $N=2$ to $672$
for $N=3$ and $3366$ for $N=4$. For uncertified runs, these counts describe processing before the common $2500$-open-node
limit, rather than the total work needed for certification; the $(50,6)$ cell is not directly comparable because it uses the larger $10000$
frontier. These results sharpen the evidence that nest count drives tree growth, consistent
with the nest dependence of the conditional complexity bound. The observed pattern does
not verify the theorem's technical assumptions or its worst-case bound, especially for the weaker PPC relaxation. It also
shows that sub-$0.3\%$ relative shifted-loss bounds remain available beyond four nests.

Appendix~\ref{app:economic-effect} reports the economic comparison of NL and calibrated MNL
pricing decisions, including the profit losses from evaluating MNL prices under NL demand.
We also compare B\&B with multistart SLSQP; Appendix~\ref{app:experimental-protocol}
gives the benchmark procedures and comparison results.
Deliberately selected challenging instances illustrate the risk of single-start local search:
SLSQP loses $1.15\%$--$48.73\%$ of profit relative to the DE incumbent. B\&B improves on
single-start SLSQP in all four cases, but none certifies its $0.5\%$ shifted-loss target.
These cases demonstrate possible local-solver failure, rather than its frequency;
Appendix~\ref{sec:local-optimum-results} gives the setup and full results.

\subsection{Product- and nest-count ablation}
\label{sec:size-ablation}
\begin{figure}[htbp]
\centering
\begin{tikzpicture}
\begin{groupplot}[
 group style={group size=2 by 2, horizontal sep=1.7cm, vertical sep=1.75cm},
 width=0.5\textwidth, height=3.7cm,
 tick label style={font=\footnotesize}, label style={font=\footnotesize},
 title style={font=\footnotesize, yshift=-3pt},
 legend style={font=\scriptsize, draw=none, fill=none, row sep=-2pt},
 legend cell align=left]
\nextgroupplot[title={Products: median wall time (s)}, xmode=log, ymode=log,
 xmin=8, xmax=620, xtick={10,25,50,100,200,500}, xticklabels={10,25,50,100,200,500},
 ymin=0.1, ymax=1000, ytick={0.1,1,10,100,1000}, yticklabels={0.1,1,10,100,1000},
 xlabel={Number of products $m$}, legend pos=north west]
\addplot[black, mark=*, line width=0.9pt, mark size=1.7pt] coordinates
 {(10,0.38) (25,0.60) (50,2.07) (100,4.47) (200,17.76) (350,65.69) (500,195.54)};
\addplot[black!55, dashed, mark=square*, line width=0.9pt, mark size=1.6pt] coordinates
 {(10,3.05) (25,6.46) (50,11.65) (100,34.98) (200,76.07) (350,181.44) (500,410.76)};
\legend{$N=2$, $N=4$}
\nextgroupplot[title={Products: median solved nodes}, xmode=log, ymode=log,
 xmin=8, xmax=620, xtick={10,25,50,100,200,500}, xticklabels={10,25,50,100,200,500},
 ymin=10, ymax=10000, ytick={10,100,1000,10000}, yticklabels={10,100,1000,10000},
 xlabel={Number of products $m$}]
\addplot[black, mark=*, line width=0.9pt, mark size=1.7pt] coordinates
 {(10,52) (25,60) (50,103) (100,55) (200,55) (350,55) (500,55)};
\addplot[black!55, dashed, mark=square*, line width=0.9pt, mark size=1.6pt] coordinates
 {(10,430) (25,553) (50,977) (100,752) (200,603) (350,562) (500,585)};
\nextgroupplot[title={Nests: median solved nodes}, ymode=log, xmin=0.5, xmax=10.5,
 xtick={1,...,10}, ymin=10, ymax=100000, ytick={10,100,1000,10000,100000},
 yticklabels={10,100,1000,$10^4$,$10^5$}, xlabel={Number of nests $N$},
 legend pos=north west]
\addplot[black, mark=*, line width=0.9pt, mark size=1.7pt] coordinates
 {(1,16) (2,81) (3,194) (4,679) (5,3239) (6,13489) (7,15012) (8,34720) (9,24665) (10,22563)};
\addplot[black!55, dashed, mark=square*, line width=0.9pt, mark size=1.6pt] coordinates
 {(1,23) (2,89) (3,303) (4,979) (5,2453) (6,10738) (7,36383) (8,23908) (9,23297) (10,25487)};
\draw[black!45, densely dotted, line width=0.6pt] (axis cs:6.5,10) -- (axis cs:6.5,100000);
\legend{$m=20$, $m=40$}
\nextgroupplot[title={Nests: median final gap (\%)}, ymode=log, xmin=0.5, xmax=10.5,
 xtick={1,...,10}, ymin=0.001, ymax=1, ytick={0.001,0.01,0.1,1},
 yticklabels={0.001,0.01,0.1,1}, xlabel={Number of nests $N$}]
\addplot[black, mark=*, line width=0.9pt, mark size=1.7pt] coordinates
 {(1,0.0072) (2,0.0080) (3,0.0090) (4,0.0099) (5,0.0100) (6,0.0100) (7,0.0305) (8,0.0364)
  (9,0.3090) (10,0.3100)};
\addplot[black!55, dashed, mark=square*, line width=0.9pt, mark size=1.6pt] coordinates
 {(1,0.0012) (2,0.0074) (3,0.0096) (4,0.0100) (5,0.0100) (6,0.0100) (7,0.0275) (8,0.1304)
  (9,0.3851) (10,0.1863)};
\draw[black!45, dashed, line width=0.5pt] (axis cs:0.5,0.01) -- (axis cs:10.5,0.01);
\draw[black!45, densely dotted, line width=0.6pt] (axis cs:6.5,0.001) -- (axis cs:6.5,1);
\end{groupplot}
\end{tikzpicture}
\caption{Effect of product and nest counts on B\&B performance. Top: product track, in which
all $28$ runs certify. Bottom: nest track; all runs certify for $N\le6$ (left of the dotted
line). Beyond it, uncertified runs end at the open-node limit; their node counts and gaps
are measured at that stop. The dashed horizontal
line marks the $0.01\%$ target. All vertical axes and the product-count axes are logarithmic.}
\label{fig:size-ablation}
\end{figure}
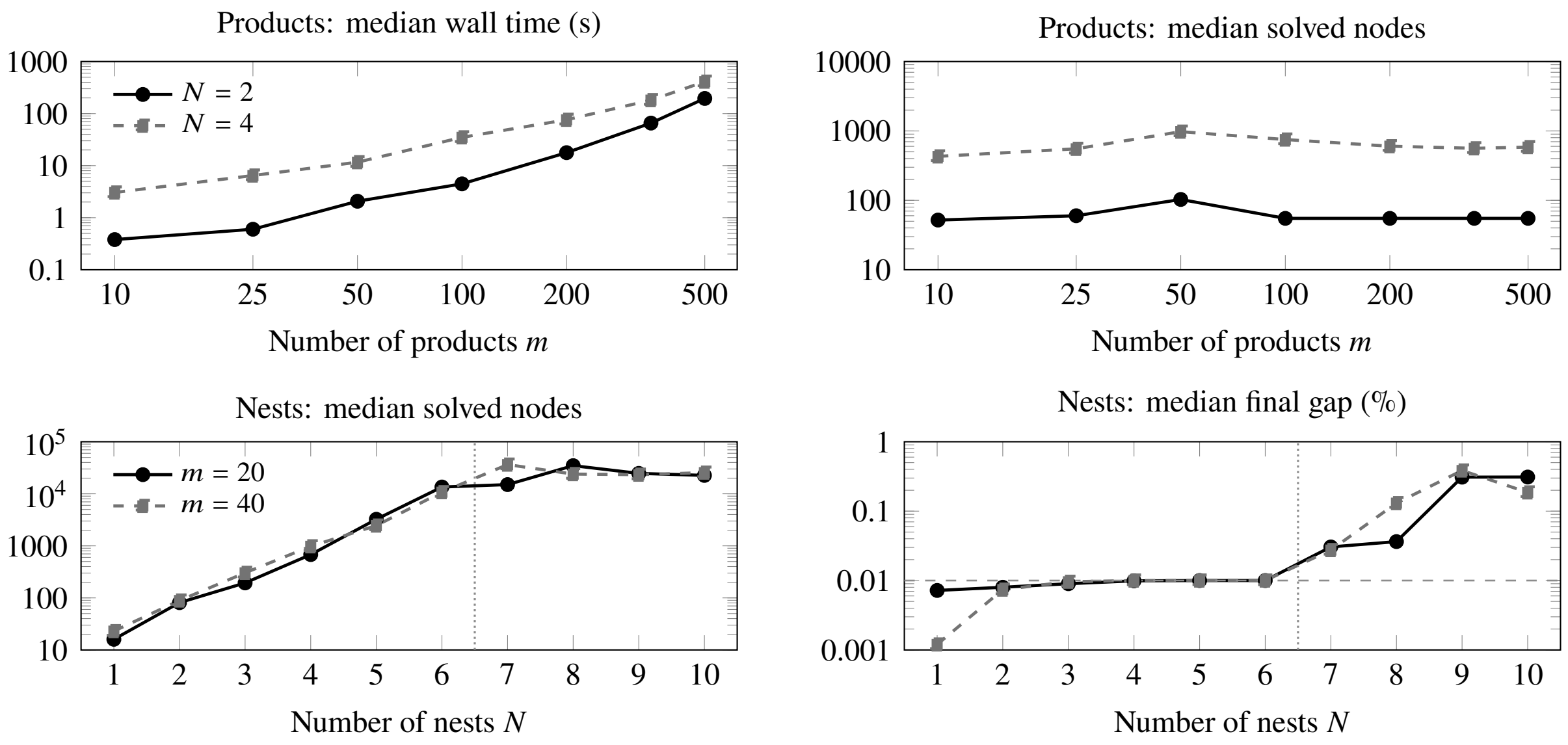

This experiment separates product-level node cost from nest-level tree
growth and tests whether the two dimensions interact. For the product track, we cross
$N\in\{2,4\}$ with $m\in\{10,25,50,100,200,350,500\}$. For the nest track, we cross
$m\in\{20,40\}$ with $N\in\{1,\ldots,10\}$. Each point has two independent replications from
the medium scale regime with LC rows. The time limit is $3600$ seconds. The open-node limit
is $2500$ in the product track and for $N\le5$; it is $10000$ for $N\ge6$ in the nest track.
Figure~\ref{fig:size-ablation} reports medians over the two replications.

Across the $68$ size-ablation runs, $55$ certify to the target gap. Every uncertified run stops
at the open-node limit, and no run approaches the time limit; the longest takes $872$ seconds.
The product track certifies all $28$ runs, up to $m=500$. Product count acts mainly through
node cost: from $m=10$ to $m=500$, median solved nodes stay between $52$ and $103$ for $N=2$
and between $430$ and $977$ for $N=4$, while median time rises from $0.38$ to $196$ seconds
and from $3.05$ to $411$ seconds, respectively. A larger fixed nest count thus raises the
tree-size constant by roughly an order of magnitude without changing the product-count trend.

The nest tracks display the sharper growth suggested by the link count in
Theorem~\ref{thm:bb-complexity}; the experiment does not establish its worst-case rate.
For both $m=20$ and $m=40$, all replications certify through $N=6$, and one of two certifies
at $N=7$; at $N=8$, one $m=20$ replication still certifies. Up to $N=6$, each added nest
multiplies the median solved-node count by a factor between $2.4$ and $5.1$: the count grows
from $16$ at $N=1$ to $13489$ at $N=6$ for $m=20$, and from $23$ to $10738$ for $m=40$,
reaching $36383$ at $(m,N)=(40,7)$. For $N\ge8$, uncertified runs in both tracks reach the $10000$-open-node
safeguard; their node counts therefore measure resource-limited searches, and the observed
final gaps need not vary monotonically. At
$N=10$, the median gaps are $0.31\%$ and $0.19\%$ for $m=20$ and $m=40$, respectively, and no
gap in either track exceeds $0.50\%$.

The nest-scale ablation is reported in Appendix~\ref{sec:sigma-ablation}.
Appendix~\ref{sec:constant-markup-experiment} compares direct constrained optimization with
projection of unconstrained CM prices.

\section{Conclusion}\label{sec:conclusion}

This paper develops an exact reformulation and a global optimization method that exploit
nest structure to return feasible constrained prices with a bound on remaining profit
improvement. The base reformulation isolates nonconvexity in $N+1$ bilinear links and represents all other
nonlinear constraints with exponential cones. Spatial B\&B then solves convex node relaxations
while retaining global bounds. Under the technical assumptions and controlled conditioning
of Theorem~\ref{thm:bb-complexity}, fixed nest count yields arithmetic work polynomial in
product count, constraint count, and inverse accuracy. A controlled positive feasible-profit
scale also gives the relative approximation scheme in Corollary~\ref{cor:fixed-nest-approx}.
The operational extensions identify which constraints preserve this nest-dependent search
dimension and which introduce additional links.

The key computational finding is that the method accommodates hundreds of products when the
number of nests is small: product count mainly increases node cost, whereas nest count drives
search-tree growth. The pricing comparisons also show the value of optimizing operational
constraints directly. Rapid search growth with nest count remains the main limitation.
Future work should strengthen nest-level relaxations and extend the approach to signed price
couplings, multiple customer segments, and dynamic inventory control.

\bibliographystyle{informs2014}
\bibliography{ref}

\clearpage
\appendix
\counterwithin{equation}{section}
\counterwithin{table}{section}
\counterwithin{figure}{section}
\counterwithin{algorithm}{section}

\pdfbookmark[1]{Appendix}{appendix-title}
\begin{center}
 {\Large\bfseries Appendix}
\end{center}
\vspace{0.5ex}

\section{Proof of Proposition~\ref{prop:reformulation-equivalence}}
\label{app:reformulation-proof}

We prove both directions of the price transformation and show why the auxiliary
inequalities do not change the optimum. For any feasible price vector, \eqref{eq:u-definition}--\eqref{eq:y-def} construct transformed
variables with
\[
 h=\alpha+\sum_{n\in\cN}z_n^{\sigma_n-1}y_n,
 \qquad k=1+\sum_{n\in\cN}z_n^{\sigma_n},
 \qquad r=h/k=L(\bp).
\]
Thus, every feasible price has an extended representation with the same shifted loss.
Conversely, let an arbitrary point of \eqref{prob:nlbc} be feasible, and recover $\bp$
by \eqref{eq:recover-price}. Its price intervals follow from \eqref{eq:u-bounds}, and its
coupling constraints follow from \eqref{con:general-price}. Define the exact nest sum
$Y_n(\bu)=\sum_{j\in\cM_n}(q_j/\beta_j)u_j\log u_j>0$. Positivity and the auxiliary
inequalities give
\[
 h\ge\alpha+\sum_n w_n
   \ge\alpha+\sum_n\rho_n y_n
   \ge\alpha+\sum_n z_n^{\sigma_n-1}Y_n(\bu),
 \qquad
 0<k\le1+\sum_n\eta_n\le1+\sum_nz_n^{\sigma_n}.
\]
Since $h\le rk$, division by the positive denominator yields
\[
 r\ge\frac{h}{k}
 \ge\frac{\alpha+\sum_n z_n^{\sigma_n-1}Y_n(\bu)}
          {1+\sum_nz_n^{\sigma_n}}
 =L(\bp).
\]
The first direction shows that the extended optimum cannot exceed the minimum shifted loss;
the second shows that it cannot fall below it. Equality of optimal values follows, and an
optimal extended point recovers an optimal price vector. Conversely, every optimal price
vector admits an optimal exact representation. Finally, $\Pi(\bp)=\alpha-L(\bp)$ gives
the stated optimal profit. \hfill$\square$

\section{Complete spatial B\&B procedure}\label{app:algorithm}

Algorithm~\ref{alg:spatial-bb} summarizes Section~\ref{sec:algorithm}. It uses node boxes
$\mathcal B_\nu$, inherited lower bounds $\underline L_\nu$, the best-bound queue $\mathcal Q$,
and the incumbent $(\widehat\bp,\overline L)$. Root processing must first obtain a feasible
incumbent and a valid lower bound; otherwise a feasibility search or an infeasible/unresolved
status is required. A failed numerical solve never proves infeasibility.
Branching uses the normalized positive residuals defined there, with envelope-error fallback
when every residual is at most $\varepsilon_{\rm bil}$. Factor widths are normalized by their
initial widths; fixed factors are excluded. Propagation discards only boxes proved infeasible.

\paragraph{Global-bound bookkeeping.}
A child inherits its parent's valid lower bound until its own relaxation is solved.
To account for tolerance-pruned boxes, let $L_{\rm disc}$ be their smallest lower bound,
initialized to $+\infty$. The global bound is
$\underline L=\min\{\overline L,L_{\rm disc},\min_{\nu\in\mathcal Q}\underline L_\nu\}$,
where $\min\emptyset=+\infty$. Exact incumbent-dominance deletions need no record, but
boxes discarded at $\overline L-\varepsilon_{\rm abs}$ may contain slightly better prices
and must be included. Tree exhaustion therefore certifies the absolute tolerance, not
necessarily zero error. Incumbent updates require original-model feasibility;
$\varepsilon_{\rm bil}$ controls branching only.

\paragraph{Persistent node model.}
One MOSEK Fusion model retains its variables, cones, sparsity, and affine expressions.
Factor bounds and envelope coefficients are parameters updated at each node
\citep[Sections~6.6 and~6.9]{mosekfusion}. Reuse saves construction work without requiring
an interior-point warm start. If numerical recovery requires reconstruction, the old model
is disposed first. Searches are serial; independent instances run in parallel.

\begin{algorithm}[!htbp]
\caption{Spatial B\&B with convex-envelope node relaxations}
\label{alg:spatial-bb}
\footnotesize
\begin{algorithmic}[1]
\Require Instance data, initial box $\mathcal B_0$, tolerances
$\varepsilon_{\rm rel},\varepsilon_{\rm abs},\varepsilon_{\rm bil}$, and resource limits
\Ensure Feasible incumbent, global lower bound, and termination status
\State Propagate initial bounds and build one persistent convex-envelope model
\State Solve the root, recover and polish an incumbent; tighten $r$ using $\overline L$ and re-solve
\State $\mathcal Q\gets\{(\mathcal B_0,\underline L_0)\}$; $L_{\rm disc}\gets+\infty$; status $\gets\textsc{resource-limit}$
\While{$\mathcal Q\ne\emptyset$ and resource limits are not reached}
  \State $\underline L\gets\min\{\overline L,L_{\rm disc},\min_{\nu\in\mathcal Q}\underline L_\nu\}$
  \If{$\overline L-\underline L\le\varepsilon_{\rm abs}$ or $\operatorname{gap}(\underline L,\overline L)\le\varepsilon_{\rm rel}$}
    \State \Return $(\widehat\bp,\overline L,\underline L,\textsc{certified})$
  \EndIf
  \State Pop the best-bound node $\nu$, propagate its box, and update the model parameters
  \State If propagation proves infeasibility, discard $\nu$ and \textbf{continue}
  \State Solve the node relaxation with numerical recovery; discard and continue only if infeasibility is proved
  \If{a valid primal--dual solution is available}
    \State Strengthen $\underline L_\nu$ with the dual bound
    \If{$\underline L_\nu\ge\overline L-\varepsilon_{\rm abs}$}
      \State $L_{\rm disc}\gets\min\{L_{\rm disc},\underline L_\nu\}$; \textbf{continue}
    \EndIf
    \State Recover prices and polish selectively; update the incumbent only after checking original feasibility
    \State Delete queued nodes whose bounds are at least the improved $\overline L$
    \State Choose the largest-residual link, or the largest envelope-error link if residuals are small
    \State Select its factor with the largest normalized interval width
  \Else
    \State Select the widest normalized factor interval; retain the inherited node bound
  \EndIf
  \If{no interval can be selected or bisected at machine precision}
    \State Reinsert $\nu$; status $\gets\textsc{spatial-resolution-limit}$; \textbf{break}
  \EndIf
  \State Bisect the interval; propagate the children and queue surviving boxes with bound $\underline L_\nu$
\EndWhile
\State If $\mathcal Q=\emptyset$, set status $\gets\textsc{tree-exhausted}$
\State $\underline L\gets\min\{\overline L,L_{\rm disc},\min_{\nu\in\mathcal Q}\underline L_\nu\}$, with $\min\emptyset=+\infty$
\State \Return $(\widehat\bp,\overline L,\underline L,\text{status})$
\end{algorithmic}
\end{algorithm}

\section{Proofs of Theorem~\ref{thm:bb-complexity} and Corollary~\ref{cor:fixed-nest-approx}}\label{app:complexity-proof}

The proof uses the standard convex-envelope and spatial-partitioning framework
\citep{mccormick1976,tawarmalani2002}. It first bounds the number of boxes needed to make
all product envelopes accurate, then uses the theorem's URR condition to convert envelope accuracy into objective accuracy. This condition is an
additional hypothesis for the present pricing formulation; envelope convergence alone does
not establish feasible-price recovery at the claimed rate. Finally, we bound the convex work
per node. For this proof, reindex the links as $s=1,\ldots,d_{\rm L}$. Let the $s$th bilinear link be $v_s=x_sy_s$, for $s=1,\ldots,d_{\rm L}$, and let its
current factor box be
$[\underline x_s,\overline x_s]\times[\underline y_s,\overline y_s]$.
The superscript $0$ denotes initial widths. Write $\Delta x_s=\overline x_s-\underline x_s$ and
$\Delta y_s=\overline y_s-\underline y_s$. The four standard envelope inequalities give the convex
hull of the graph of $x_sy_s$ over this rectangle. The maximum vertical distance between the
product graph and either envelope surface is bounded by
\begin{equation}
 \frac{\Delta x_s\Delta y_s}{4}.
 \label{eq:mccormick-error-bound}
\end{equation}
This follows by translating and scaling the rectangle to $[0,1]^2$: there, the largest gap
between $xy$ and its convex or concave envelope is $1/4$, attained at the center. Scaling
back multiplies the gap by $\Delta x_s\Delta y_s$.
This also explains why each bilinear link contributes one, rather than two, to the final
accuracy exponent. Reducing both factor widths to order $\sqrt{\theta}$ makes their
product error order $\theta$. Each factor then needs order $\theta^{-1/2}$ intervals,
giving order $\theta^{-1}$ rectangles per pair of factors.

Fix a desired product tolerance $\theta>0$. For each link, let
\[
 t_s=\max\left\{1,\sqrt{\frac{\Delta x_s^0\Delta y_s^0}{4\theta}}\right\},
 \qquad r_s=2^{\lceil\log_2 t_s\rceil}.
\]
Bisect each initial factor interval $\log_2r_s$ times. Both intervals then have $r_s$
equal pieces, and each resulting rectangle satisfies
$\Delta x_s\Delta y_s/4\le\theta$. The power-of-two choice makes this partition
realizable by the binary refinement assumed in the theorem. Since $r_s<2t_s$, the number
of leaf boxes across all $d_{\rm L}$ factor pairs is at most
\begin{equation}
 \prod_{s=1}^{d_{\rm L}}r_s^2
 \le4^{d_{\rm L}}\prod_{s=1}^{d_{\rm L}}
       \max\left\{1,\frac{\Delta x_s^0\Delta y_s^0}{4\theta}\right\}
 =O(\theta^{-d_{\rm L}}).
 \label{eq:partition-cardinality}
\end{equation}
If a factor occurs in several links, retaining only mutually consistent interval assignments
can only reduce this count. A full binary tree has fewer than twice as many total nodes as
leaves, giving the same $O(\theta^{-d_{\rm L}})$ bound for all processed nodes.

By the assumed URR condition, choosing
$\theta=\varepsilon/(2\kappa d_{\rm L})$ makes each surviving leaf either certifiably
infeasible or supplies a feasible point within $\varepsilon/2$ of its exact relaxation
optimum. Reserve the remaining $\varepsilon/2$ for the difference between that optimum
and the computed valid dual lower bound. Each recovered point is then within
$\varepsilon$ of its leaf bound. Taking the point recovered from a leaf with the smallest
valid bound gives an incumbent within $\varepsilon$ of that global bound. Substituting this
value of $\theta$ into
\eqref{eq:partition-cardinality}, collecting the initial widths and $2\kappa d_{\rm L}$ into
$C_{\rm B}$, and using $d_{\rm L}=N+1$ proves \eqref{eq:bb-node-complexity}.

It remains to bound the work at one node. After convex-envelope replacement, all constraints
are affine or log-domain convex, so the node problem is convex. The standard logarithmically
homogeneous barrier for one three-dimensional exponential cone has constant barrier
parameter three \citep[Section~2]{dahlandersen2022}. Adding the barriers for the $n_e$ exponential cones and the scalar conic
components gives total parameter $\vartheta=O(n_v+n_B+3n_e)$. Under the stated regularity assumptions, self-concordant-barrier theory
\citep{nesterovnemirovskii1994} gives
$O(\bigl.\sqrt{\vartheta}\log(\vartheta/\tau)\bigr)$ Newton steps to reach primal--dual accuracy $\tau$.
Using dense linear algebra, forming and factorizing the Newton system costs $O(n_K^3)$
arithmetic operations per step. Setting $\tau=\Theta(\varepsilon)$, adding the auxiliary work
$T_{\rm aux}(\varepsilon)$, and multiplying the per-node cost by
\eqref{eq:bb-node-complexity} yields
\eqref{eq:bb-total-complexity}. Sparse factorization and persistent model reuse may improve the
practical cost but do not change this worst-case conclusion. The displayed node count is
polynomial in $1/\varepsilon$ for fixed $d_{\rm L}$ and fixed $C_{\rm B}$; an input-dimension polynomial
claim additionally requires the widths and $\kappa$ collected in $C_{\rm B}$ to be polynomially
bounded. \hfill$\square$

\paragraph{Proof of Corollary~\ref{cor:fixed-nest-approx}.}
Apply Theorem~\ref{thm:bb-complexity} with absolute tolerance $\varepsilon=\gamma R_0$.
Since $R_0\le\Pi^*$ and shifted-loss and profit errors coincide,
$\Pi^*-\Pi(\widehat\bp)\le\gamma R_0\le\gamma\Pi^*$. The node factor becomes
$(C_{\rm B}/(\gamma R_0))^{N+1}$, which is polynomial under the stated family-level
bounds and fixed $N$. \hfill$\square$

\section{Detailed computational results}\label{app:detailed-results}

Tables~\ref{tab:final-family-performance}--\ref{tab:ec-size-ablation} retain the constraint-system
comparison and product- and nest-count results supporting Section~\ref{sec:experiments}.
Times are in seconds and gaps are shifted-loss percentages \eqref{eq:bb-relative-gap}.
Results use two decimal places; $0.00$ can denote a small positive value. Certification
uses unrounded values. All $160$ matched runs certify, as do all product-count cases through
$500$ products with two or four nests. Certification becomes harder as nest count increases;
uncertified runs in these tables stop at the open-node limit and retain valid bounds.
The economic comparisons and other ablations appear with their discussions in the
following appendices, avoiding repeated tables here.

\begin{table}[htbp]
\centering
\caption{B\&B median wall time (seconds) and median solved nodes in the matched campaign; all
$160$ runs certify the $0.01\%$ target. LC: linear price constraints; PPC:
purchase-probability constraints; LC+PPC: their combination; POC: pairwise price-order
constraints. Each NL cell pools the two scale regimes and two replications (one for POC);
each MNL row pools the matched $N=2$ and $N=4$ partition labels. POC was tested only at
$m=5,10$; its dashes at $m=20$ indicate untested cases.}
\label{tab:final-family-performance}
\footnotesize
\setlength{\tabcolsep}{4.2pt}
\begin{tabular}{@{}lrrrrrrrrrr@{}}
\toprule
& & & \multicolumn{2}{c}{LC} & \multicolumn{2}{c}{PPC} & \multicolumn{2}{c}{LC+PPC}
& \multicolumn{2}{c}{POC}\\
\cmidrule(lr){4-5}\cmidrule(lr){6-7}\cmidrule(lr){8-9}\cmidrule(l){10-11}
Model & $m$ & $N$ & Time & Nodes & Time & Nodes & Time & Nodes & Time & Nodes\\
\midrule
NL & 5 & 2 & 0.44 & 89 & 0.41 & 82 & 0.49 & 90 & 0.15 & 42\\
NL & 5 & 4 & 3.03 & 455 & 3.87 & 618 & 4.35 & 499 & 0.71 & 158\\
NL & 10 & 2 & 0.46 & 55 & 0.44 & 43 & 0.60 & 64 & 0.44 & 88\\
NL & 10 & 4 & 6.82 & 637 & 5.67 & 636 & 6.72 & 632 & 2.15 & 421\\
NL & 20 & 2 & 0.86 & 91 & 0.85 & 60 & 0.94 & 75 & -- & --\\
NL & 20 & 4 & 14.49 & 1557 & 29.59 & 3922 & 17.81 & 1819 & -- & --\\
\addlinespace
MNL & 5 & -- & 0.08 & 2 & 0.10 & 2 & 0.10 & 2 & 0.06 & 2\\
MNL & 10 & -- & 0.17 & 2 & 0.28 & 2 & 0.28 & 2 & 0.13 & 2\\
MNL & 20 & -- & 0.36 & 2 & 0.56 & 2 & 0.64 & 2 & -- & --\\
\bottomrule
\end{tabular}
\end{table}

\begin{table}[htbp]
\centering
\caption{B\&B performance across problem sizes and constraint systems in Section~\ref{sec:final-results}. ``Cert.'' gives
certificates over the number of runs. Parentheses contain the maximum gap, mean time, and mean
node count. Each NL row contains the four constraint systems; each MNL row pools the systems
over matched partition labels. The $(m,N)=(50,6)$ cell uses the $10000$-open-node limit.}
\label{tab:ec-large-instance}
\footnotesize
\setlength{\tabcolsep}{5pt}
\begin{tabular}{@{}lrrrrrrr@{}}
\toprule
Model & $m$ & $N$ & Runs & Cert. & \shortstack{Gap (\%)\\med. (max.)} &
\shortstack{Time (s)\\med. (mean)} & \shortstack{Nodes\\med. (mean)}\\
\midrule
NL & 50 & 2 & 4 & 4/4 & 0.01 (0.01) & 1.90 (2.45) & 83 (84) \\
NL & 50 & 3 & 4 & 4/4 & 0.01 (0.01) & 12.50 (12.08) & 746 (658) \\
NL & 50 & 4 & 4 & 1/4 & 0.06 (0.08) & 68.75 (59.69) & 6316 (4948) \\
NL & 50 & 5 & 4 & 1/4 & 0.10 (0.16) & 89.95 (76.55) & 5643 (4792) \\
NL & 50 & 6 & 4 & 0/4 & 0.18 (0.28) & 377.25 (367.90) & 22563 (22519) \\
NL & 100 & 2 & 4 & 4/4 & 0.01 (0.01) & 5.90 (6.59) & 76 (76) \\
NL & 100 & 3 & 4 & 4/4 & 0.01 (0.01) & 24.50 (25.11) & 672 (560) \\
NL & 100 & 4 & 4 & 4/4 & 0.01 (0.01) & 112.08 (103.00) & 3366 (2930) \\
NL & 100 & 5 & 4 & 1/4 & 0.03 (0.07) & 180.18 (187.00) & 5663 (4362) \\
NL & 100 & 6 & 4 & 0/4 & 0.09 (0.13) & 269.97 (304.20) & 6018 (5991) \\
NL & 150 & 2 & 4 & 4/4 & 0.01 (0.01) & 11.64 (16.26) & 84 (80) \\
NL & 150 & 3 & 4 & 4/4 & 0.01 (0.01) & 48.17 (43.00) & 698 (618) \\
NL & 150 & 4 & 4 & 1/4 & 0.03 (0.04) & 371.46 (304.60) & 5993 (4943) \\
NL & 150 & 5 & 4 & 1/4 & 0.03 (0.05) & 345.42 (371.90) & 5756 (4462) \\
NL & 150 & 6 & 4 & 0/4 & 0.10 (0.21) & 354.43 (385.10) & 5624 (5773) \\
\addlinespace
MNL & 50 & -- & 12 & 12/12 & 0.00 (0.00) & 1.31 (1.60) & 2 (2) \\
MNL & 100 & -- & 12 & 12/12 & 0.00 (0.00) & 4.76 (5.31) & 2 (2) \\
MNL & 150 & -- & 12 & 12/12 & 0.00 (0.00) & 8.64 (11.16) & 2 (2) \\
\bottomrule
\end{tabular}
\end{table}

\begin{table}[htbp]
\centering
\caption{B\&B product- and nest-count ablation (Section~\ref{sec:size-ablation}). Gap entries are median
(maximum) percentages; time and solved-node entries are medians over two replications.}
\label{tab:ec-size-ablation}
\footnotesize
\setlength{\tabcolsep}{4.5pt}
\begin{tabular}{@{}rrrrrrrrr@{}}
\toprule
& \multicolumn{4}{c}{$N=2$} & \multicolumn{4}{c}{$N=4$}\\
\cmidrule(lr){2-5}\cmidrule(l){6-9}
$m$ & Cert. & Gap & Time & Solved & Cert. & Gap & Time & Solved\\
\midrule
10 & 2/2 & 0.01 (0.01) & 0.38 & 52 & 2/2 & 0.01 (0.01) & 3.05 & 430 \\
25 & 2/2 & 0.01 (0.01) & 0.60 & 60 & 2/2 & 0.01 (0.01) & 6.46 & 553 \\
50 & 2/2 & 0.01 (0.01) & 2.07 & 103 & 2/2 & 0.01 (0.01) & 11.65 & 977 \\
100 & 2/2 & 0.01 (0.01) & 4.47 & 55 & 2/2 & 0.01 (0.01) & 34.98 & 752 \\
200 & 2/2 & 0.01 (0.01) & 17.76 & 55 & 2/2 & 0.01 (0.01) & 76.07 & 603 \\
350 & 2/2 & 0.01 (0.01) & 65.69 & 55 & 2/2 & 0.01 (0.01) & 181.44 & 562 \\
500 & 2/2 & 0.01 (0.01) & 195.54 & 55 & 2/2 & 0.01 (0.01) & 410.76 & 585 \\
\addlinespace
& \multicolumn{4}{c}{$m=20$} & \multicolumn{4}{c}{$m=40$}\\
\cmidrule(lr){2-5}\cmidrule(l){6-9}
$N$ & Cert. & Gap & Time & Solved & Cert. & Gap & Time & Solved\\
\midrule
1 & 2/2 & 0.01 (0.01) & 0.23 & 16 & 2/2 & 0.00 (0.00) & 0.49 & 23 \\
2 & 2/2 & 0.01 (0.01) & 0.92 & 81 & 2/2 & 0.01 (0.01) & 0.99 & 89 \\
3 & 2/2 & 0.01 (0.01) & 1.92 & 194 & 2/2 & 0.01 (0.01) & 4.39 & 303 \\
4 & 2/2 & 0.01 (0.01) & 3.95 & 679 & 2/2 & 0.01 (0.01) & 16.17 & 979 \\
5 & 2/2 & 0.01 (0.01) & 20.78 & 3239 & 2/2 & 0.01 (0.01) & 44.49 & 2453 \\
6 & 2/2 & 0.01 (0.01) & 163.30 & 13489 & 2/2 & 0.01 (0.01) & 216.20 & 10738 \\
7 & 1/2 & 0.03 (0.05) & 141.80 & 15012 & 1/2 & 0.03 (0.05) & 648.40 & 36383 \\
8 & 1/2 & 0.04 (0.06) & 373.00 & 34720 & 0/2 & 0.13 (0.19) & 339.10 & 23908 \\
9 & 0/2 & 0.31 (0.48) & 288.80 & 24665 & 0/2 & 0.39 (0.50) & 326.20 & 23297 \\
10 & 0/2 & 0.31 (0.37) & 271.20 & 22563 & 0/2 & 0.19 (0.24) & 393.90 & 25487 \\
\bottomrule
\end{tabular}
\end{table}

\section{Economic effect of NL versus MNL}\label{app:economic-effect}

This section examines the economic effect of replacing NL demand by a calibrated MNL
approximation, using the matched instances of Section~\ref{sec:final-results}. All $80$ matched NL--MNL pairs return feasible incumbents. Write $\Pi_{\rm NL}(\bp)$ and $\Pi_{\rm MNL}(\bp)$ for profit per customer under the
two demand models, and $\widehat\bp^{\rm NL}$ and $\widehat\bp^{\rm MNL}$ for their
certified incumbent price vectors. The policy loss evaluates both prices under the same NL
demand and equals
$100[\Pi_{\rm NL}(\widehat\bp^{\rm NL})-\Pi_{\rm NL}(\widehat\bp^{\rm MNL})]/
\Pi_{\rm NL}(\widehat\bp^{\rm NL})$. For each pair,
Table~\ref{tab:nl-mnl-economic} summarizes the pairwise percentage differences between NL and MNL
incumbent profits and the policy losses under NL.
Each profit difference is relative to the calibrated-MNL incumbent profit,
$100[\Pi_{\rm NL}(\widehat\bp^{\rm NL})-\Pi_{\rm MNL}(\widehat\bp^{\rm MNL})]/
\Pi_{\rm MNL}(\widehat\bp^{\rm MNL})$. Every MNL price
vector remains feasible in its matched NL instance.
This is an instance-specific outcome: the unchanged price bounds and price rules transfer
directly, but PPC must be checked under NL demand. Matching
conversion at midpoint prices alone does not guarantee feasibility at other prices.

\begin{table}[htbp]
\centering
\caption{Economic effect of NL demand relative to calibrated MNL. Constraint-system
abbreviations follow Table~\ref{tab:final-family-performance}.}
\label{tab:nl-mnl-economic}
\footnotesize
\begin{tabular}{@{}lrrrrr@{}}
\toprule
System & Pairs & Mean NL--MNL & Median NL--MNL & Mean policy loss & Max. policy loss\\
       &       & profit (\%) & profit (\%)   & under NL (\%)   & under NL (\%)\\
\midrule
LC     & 24 & 10.16 & 4.79  & 2.44  & 9.56\\
PPC    & 24 & 17.45 & 14.50 & 10.18 & 34.42\\
LC+PPC & 24 & 11.32 & 4.04  & 5.81  & 28.49\\
POC    &  8 & 20.42 & 20.77 & 6.06  & 17.09\\
\midrule
Overall & 80 & 13.72 & 6.49 & 6.14 & 34.42\\
\bottomrule
\end{tabular}
\end{table}

Under the same NL demand, MNL prices lose $6.14\%$ of profit on average (median $2.02\%$,
maximum $34.42\%$), with the largest mean loss under PPC. Separately, optimized NL profit
exceeds calibrated-MNL profit in every pair; this cross-model difference is not a pricing
policy gain. NL prices average $0.199$ units higher. The small replication counts support
controlled comparisons rather than population estimates.

\section{Nest-scale ablation}
\label{sec:sigma-ablation}

We fix $m=50$, $N=2$, and cross mean nest scale ($0.25$ or $0.75$) with scale
separation ($0.10$ or $0.30$), giving the four scenarios in
Table~\ref{tab:sigma-ablation}. Each of three replications shares products, nests, costs,
within-nest sensitivities, bounds, and LC rows across scenarios. Intercept shifts hold
midpoint conversion at $0.50$, but do not equalize individual choice probabilities.
Scales vary at fixed $\beta_j$, not fixed conventional utility $V_j$
(Section~\ref{sec:model}). The time limit is $600$ seconds.

\begin{table}[htbp]
\centering
\caption{Nest-scale ablation. Entries are medians over three replications, with ranges in
parentheses.}
\label{tab:sigma-ablation}
\footnotesize
\setlength{\tabcolsep}{3pt}
\begin{tabular}{@{}lcrrrrr@{}}
\toprule
Scenario & $(\sigma_1,\sigma_2)$ & Cert. & \shortstack{Time (s)\\med. (range)} &
\shortstack{Gap (\%)\\med. (range)} & \shortstack{Solved\\nodes} & \shortstack{Open\\nodes}\\
\midrule
Small--close & $(0.20,0.30)$ & 3/3 & 3.84 (3.80--4.10) & 0.01 (0.01--0.01) & 192 & 11\\
Small--different & $(0.10,0.40)$ & 3/3 & 6.17 (4.79--6.67) & 0.01 (0.01--0.01) & 208 & 19\\
Large--close & $(0.70,0.80)$ & 3/3 & 1.99 (1.97--2.10) & 0.01 (0.01--0.01) & 72 & 5\\
Large--different & $(0.60,0.90)$ & 3/3 & 2.02 (1.90--2.06) & 0.01 (0.00--0.01) & 84 & 7\\
\bottomrule
\end{tabular}
\end{table}

All twelve cases certify in under seven seconds. Small scales require roughly two to
three times as many nodes as large scales, consistent with weaker relaxations; the effect
is modest with two nests and bound propagation. Greater scale separation slightly enlarges
both the solved-node count and the open-node frontier.

\section{Value of direct constrained optimization}
\label{sec:constant-markup-experiment}

We test whether projecting unconstrained CM prices can replace direct
constrained optimization. A CM $t$ sets $p_j=c_j+t$ for every product
\citep{lihuh2011,gallegowang2014}. The $16$ instances cross $m\in\{20,50\}$,
$\sigma_n=\sigma\in\{0.375,0.875\}$, and LC, PPC, LC+PPC, and POC, with $N=2$ and
$\beta_j=1$. Bounds $[c_j,c_j+6]$ admit the unconstrained solution, but each operational
system excludes it.

We compare its minimum-distance feasible projection (CM-Proj), the highest-profit feasible
common markup (FCM), SLSQP \citep{kraft1994} initialized at CM-Proj (CM-Proj+SLSQP), and B\&B.
CM-Proj minimizes the Euclidean distance $\|\bp-\bp^{\rm CM}\|_2$ to the unconstrained CM price
vector $\bp^{\rm CM}$ over the feasible set, including the price bounds. SLSQP solves this
projection from four starts (the CM vector clipped to the price box, the box midpoint, and the
lower and upper price bounds), and the closest feasible result is kept. For LC and POC the
projection is a convex quadratic program; for PPC and LC+PPC, the lower PPC
bound can make the feasible set nonconvex in prices, so CM-Proj is the closest local projection
found. The common
budget is $600$ seconds and B\&B targets $10^{-4}$. A method $b$'s profit loss is
$100(\Pi^{\rm BB}-\Pi^b)/\Pi^{\rm BB}$, where $\Pi^{\rm BB}$ and $\Pi^b$ are feasible incumbent profits;
this comparison is distinct from the certified shifted-loss gap.

\begin{figure}[htbp]
\centering
\begin{tikzpicture}
\begin{groupplot}[
 group style={group size=2 by 2, horizontal sep=1.1cm, vertical sep=1.75cm},
 width=0.5\textwidth, height=3.5cm,
 symbolic x coords={LC,PPC,LC+PPC,POC}, xtick=data, enlarge x limits=0.18, ymin=0,
 tick label style={font=\footnotesize}, label style={font=\footnotesize},
 title style={font=\footnotesize, yshift=-3pt},
 legend style={font=\scriptsize, draw=none, fill=none, row sep=-2pt},
 legend cell align=left]
\nextgroupplot[title={$\sigma=0.375$}, ylabel={Profit loss (\%)}, ybar, bar width=7pt,
 ymax=35, ytick={0,10,20,30}, error bars/y dir=plus, error bars/y explicit,
 error bars/error mark options={rotate=90, mark size=2pt}]
\addplot[draw=black, fill=black!65] coordinates
 {(LC,13.71) +- (0,4.07) (PPC,26.29) +- (0,6.32) (LC+PPC,18.28) +- (0,1.05) (POC,0.05) +- (0,0.02)};
\addplot[draw=black, fill=white] coordinates
 {(LC,7.77) +- (0,1.47) (PPC,0.00) +- (0,0.00) (LC+PPC,6.34) +- (0,1.10)};
\node[font=\scriptsize, anchor=south] at ([xshift=4.5pt]axis cs:POC,0) {n/a};
\nextgroupplot[title={$\sigma=0.875$}, ybar, bar width=7pt, ymax=35, ytick={0,10,20,30},
 yticklabels={}, legend pos=north east, area legend, error bars/y dir=plus, error bars/y explicit,
 error bars/error mark options={rotate=90, mark size=2pt}]
\addplot[draw=black, fill=black!65] coordinates
 {(LC,3.11) +- (0,1.74) (PPC,7.08) +- (0,0.80) (LC+PPC,3.13) +- (0,0.54) (POC,0.24) +- (0,0.06)};
\addplot[draw=black, fill=white] coordinates
 {(LC,13.58) +- (0,1.64) (PPC,0.02) +- (0,0.01) (LC+PPC,13.06) +- (0,1.91)};
\node[font=\scriptsize, anchor=south] at ([xshift=4.5pt]axis cs:POC,0) {n/a};
\legend{CM-Proj, FCM}
\nextgroupplot[title={Median wall time (s)}, ybar, bar width=7pt, ymax=20,
 ytick={0,5,10,15,20}, legend pos=north west, area legend]
\addplot[draw=black, fill=black!35] coordinates
 {(LC,5.48) (PPC,3.65) (LC+PPC,7.42) (POC,18.72)};
\addplot[draw=black, pattern=north east lines] coordinates
 {(LC,3.21) (PPC,2.10) (LC+PPC,3.55) (POC,4.47)};
\legend{$\sigma=0.375$, $\sigma=0.875$}
\nextgroupplot[title={Median solved nodes}, ybar, bar width=7pt, ymax=1800,
 ytick={0,500,1000,1500}, yticklabels={0,500,1000,1500}]
\addplot[draw=black, fill=black!35] coordinates
 {(LC,389) (PPC,319) (LC+PPC,458) (POC,1597)};
\addplot[draw=black, pattern=north east lines] coordinates
 {(LC,225) (PPC,126) (LC+PPC,170) (POC,284)};
\end{groupplot}
\node[anchor=south, font=\footnotesize\bfseries] at ([xshift=0.55cm, yshift=12pt]group c1r1.north east)
 {(a) Feasible-solution quality};
\node[anchor=south, font=\footnotesize\bfseries] at ([xshift=0.55cm, yshift=12pt]group c1r2.north east)
 {(b) B\&B bounding performance};
\end{tikzpicture}
\caption{CM benchmarks versus constrained optimization with $N=2$; medians over
$m\in\{20,50\}$. (a) Profit loss relative to the B\&B incumbent; whiskers reach the maximum,
and n/a marks a missing FCM. CM-Proj+SLSQP, with loss at most
$1.6\times10^{-7}\%$, is omitted. (b) All $16$ B\&B runs certify the $0.01\%$ target.}
\label{fig:constant-markup-comparison}
\end{figure}
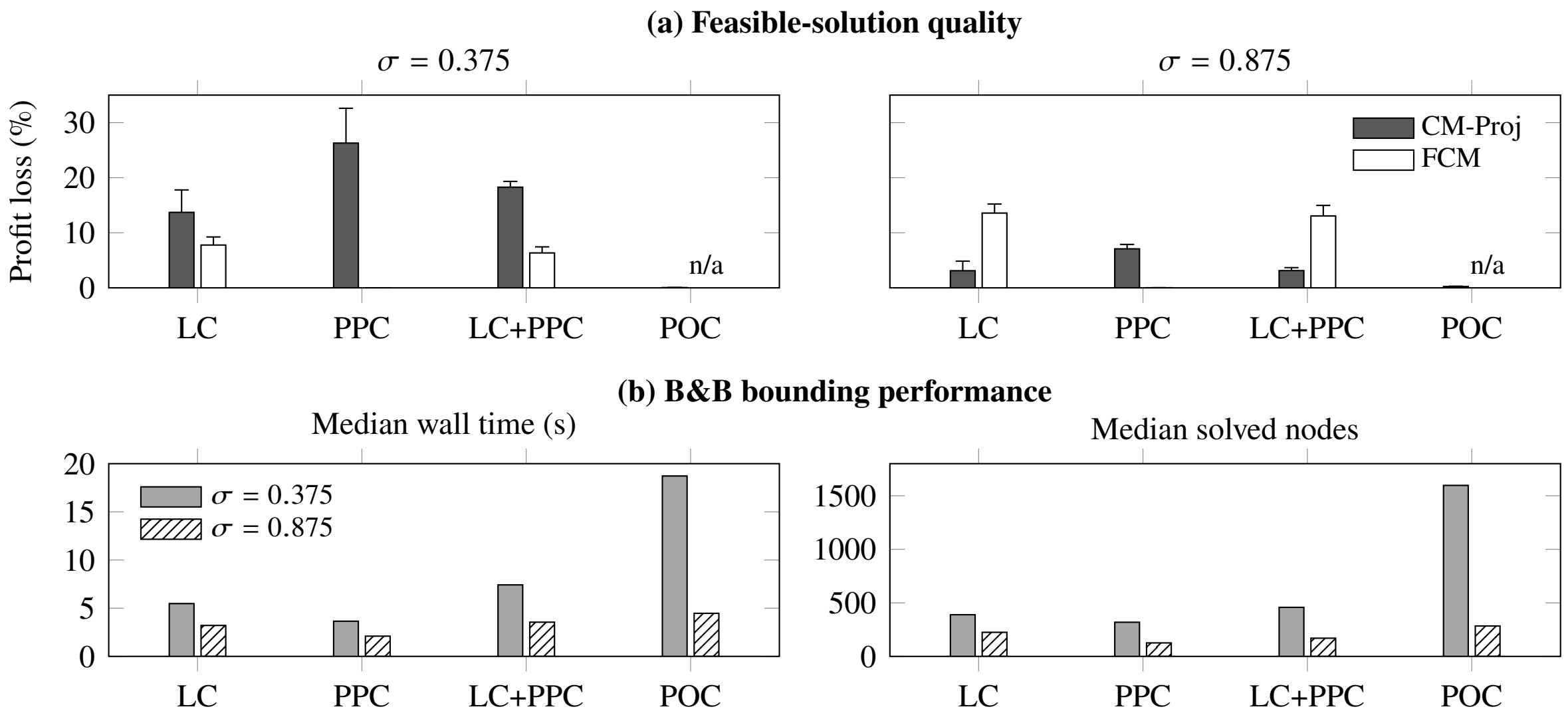

Projection is fast (median $0.16$ seconds), but CM-Proj loses $5.56\%$ of profit at the
median and up to $32.60\%$. FCM is feasible in twelve cases, with median loss $6.86\%$
and maximum $15.22\%$; no POC case admits an FCM because the required orders
reverse selected cost rankings. Strong correlation increases projection losses.
CM-Proj+SLSQP nearly matches B\&B (maximum loss $1.6\times10^{-7}\%$, median runtime
$6.57$ seconds), but provides no global bound. B\&B certifies all sixteen cases within
budget. Thus, CM prices can initialize a search, while direct constraint
handling delivers the gains and B\&B supplies the certificate.

\section{Detailed experimental protocol}\label{app:experimental-protocol}

This section supplements Section~\ref{sec:experimental-design}; departures from the default
generator and budgets are stated with each experiment.

\paragraph{Instance generation.}
We write $\mathrm{Uniform}[a,b]$ for a uniform distribution and $\mathrm{Normal}(\mu,s^2)$ for a
normal distribution with mean $\mu$ and variance $s^2$.
Draw $c_j\sim\mathrm{Uniform}[0.5,1.5]$, set $\underline p_j=c_j$, and draw
$\overline p_j-c_j\sim\mathrm{Uniform}[1,3]$. Set
$\alpha=\max_j(\overline p_j-c_j)+0.5$. Randomly permute products into approximately equal
nests, draw $\log\beta_j\sim\mathrm{Normal}(0,0.25^2)$ and raw intercepts from
$\mathrm{Normal}(0,0.5^2)$, and draw nest scales uniformly from the regimes in
Section~\ref{sec:experimental-design}. Shift intercepts to obtain midpoint conversion $0.50$.
For LC, nonzero $B$ entries are uniform on $[0,1]$ and $\bd=B\bp^{\rm mid}$; PPC uses a
conversion band of width $0.10$ centered there. POC uses one to three price-order rules and
common sensitivity uniform on $[0.75,1.25]$. The implemented PPC relaxation bounds $k$ by
\eqref{eq:total-share-k-bounds} without \eqref{con:probability-power-link}. Every feasible
price admits $k=D$, preserving a valid lower bound; recovery checks the original band.
This weaker relaxation does not inherit the exact-band complexity bound.

\paragraph{Matched campaign and benchmarks.}
The $144$ LC/PPC/LC+PPC runs cross $m\in\{5,10,20\}$, $N\in\{2,4\}$, strong-correlation
and near-MNL scales, two replications, and NL/MNL demand. Another $16$ POC runs use
$m\in\{5,10\}$ and one replication. MNL sets $\sigma_n=1$ and recalibrates intercepts to
match midpoint conversion, keeping other inputs fixed. Multistart SLSQP
\citep{kraft1994} retains the best feasible result from several initial prices; DE
\citep{stornprice1997} uses local polishing. SLSQP matches all $160$ B\&B incumbents
(maximum relative shifted-loss difference $4.0\times10^{-9}\%$); DE's maximum difference
across $112$ comparisons is $5.0\times10^{-4}\%$. Median times are $3.06$ and $9.63$ seconds,
respectively. Neither benchmark supplies a global bound.

\paragraph{Settings and verification.}
In addition to Section~\ref{sec:experimental-design}'s defaults, the bilinear-residual threshold
is $10^{-5}$ and the open-node limit rises to $10000$ when $N\ge6$ and $m\le50$, unless
specified otherwise. The software is Python~3.13.9, SciPy~1.16.3 \citep{virtanen2020}, and
MOSEK~11.2.4 \citep{mosekfusion}, on Windows~11 (Intel Core i7-14700, 20 cores, 64~GB), with
up to twelve concurrent workers. The study reports that code, generators, and raw results
are maintained with the manuscript. Its reported audit regenerates all $356$ instances,
checks incumbent feasibility, and recomputes objectives and gaps; all pass, with maximum
violation $4.0\times10^{-12}$. It does not independently verify conic dual bounds or use
outward-rounded interval arithmetic. Status \texttt{optimal\_gap} and tree exhaustion count as
certified to the stopping tolerance. Uncertified runs stop at the open-node limit, except
one time-limited run on a challenging local-optimization instance.

\section{Challenging instances for local optimization}
\label{sec:local-optimum-results}

This study uses deliberately selected challenging instances to examine whether a standard
single-start constrained solver can return economically inferior feasible prices.
It demonstrates possibility, not frequency in a representative population.

These cases use much wider price intervals and more heterogeneous sensitivities than the
main campaign. The four instances have $m=8$, $N=3$, nest cardinalities $(3,3,2)$,
$c_j\sim\mathrm{Uniform}[0,2]$, price floors at cost, markup-interval widths
$\Omega_j\sim\mathrm{Uniform}[10,30]$, attraction intercepts sampled from $\mathrm{Uniform}[-8,8]$,
$\log\beta_j\sim\mathrm{Normal}(0,1.5^2)$, and $\sigma_n\sim\mathrm{Uniform}[0.05,0.50]$. Scanning
seeds upward from a fixed value, we retain the first four instances for which SLSQP from one
fixed start ends at least $0.5\%$ below the best of $20$ SLSQP starts. We compare B\&B
implemented with MOSEK \citep{mosekfusion}, SLSQP from the fixed start \citep{kraft1994},
and DE with local polishing. B\&B uses a $0.5\%$ target and a
$600$-second limit.

\begin{table}[htbp]
\centering
\caption{Profit on the selected challenging instances for local optimization ($m=8$, $N=3$). B\&B
reports its incumbent and final shifted-loss gap at termination; SLSQP loss is measured
relative to the DE incumbent.}
\label{tab:local-optimum-profit}
\footnotesize
\setlength{\tabcolsep}{4pt}
\begin{tabular}{@{}ccrrrrr@{}}
\toprule
$[\min_j\beta_j,\max_j\beta_j]$ & $[\min_n\sigma_n,\max_n\sigma_n]$ & SLSQP & DE & B\&B &
\shortstack{SLSQP loss\\vs. DE (\%)} & \shortstack{B\&B\\gap (\%)}\\
\midrule
$[0.089,3.399]$ & $[0.362,0.451]$ & 0.95 & 1.11 & 1.11 & 14.52 & 30.27\\
$[0.656,3.771]$ & $[0.331,0.450]$ & 1.71 & 1.79 & 1.79 & 4.46 & 43.52\\
$[0.605,4.070]$ & $[0.223,0.470]$ & 0.89 & 0.90 & 0.89 & 1.15 & 25.38\\
$[0.746,6.803]$ & $[0.306,0.353]$ & 1.49 & 2.90 & 2.90 & 48.73 & 63.03\\
\bottomrule
\end{tabular}
\end{table}

Single-start SLSQP loses $1.15\%$--$48.73\%$ of profit relative to DE on these selected
cases. B\&B improves every single-start incumbent using unrounded profits; its differences
from DE on cases two and three are $0.02\%$ and $0.77\%$, which table rounding can obscure.
No B\&B run certifies the $0.5\%$ target: three hit the open-node limit after $455$--$598$
seconds and one hits the time limit, with $9854$--$16310$ solved nodes and shifted-loss
gaps of $25.4\%$--$63.0\%$. The experiment illustrates inferior local solutions and the
benefit of relaxation-guided polishing, not their population frequency or successful
certification on wide, poorly scaled boxes.

\end{document}